\documentclass[a4paper,twoside]{article}

\usepackage{amsthm,amssymb,amsmath,mathscinet,booktabs}
\usepackage[bookmarks=true,bookmarksopen=true,colorlinks=true,citecolor=blue,linkcolor=blue,urlcolor=blue]{hyperref}
\usepackage[dvipsnames]{xcolor}
\usepackage{xspace}

\numberwithin{table}{section}
\numberwithin{equation}{section}
\theoremstyle{plain}
\newtheorem{theorem}{Theorem}[section]
\newtheorem{lemma}[theorem]{Lemma}
\newtheorem{corollary}[theorem]{Corollary}

\theoremstyle{definition}

\theoremstyle{remark}
\newtheorem{remark}[theorem]{Remark}

\newtheorem{assumption}[theorem]{Assumption}

\newcommand{\dx}{\,\mathrm{d}x}

\newcommand{\Pb}{\mbox{\rm (P)}\xspace}
\newcommand{\Pbh}{\mbox{\rm (P$_h$)}\xspace}
\newcommand{\uad}{U_{\rm ad}}
\newcommand{\Uh}{U_h}
\newcommand{\uadh}{U_{h,\rm ad}}
\newcommand{\proj}{\operatorname{Proj}}

\newcommand{\umin}{u_{a}}
\newcommand{\umax}{u_{b}}
\newcommand{\ps}{{\hat p}}
\newcommand{\qs}{{\hat q}}
\newcommand{\xS}{S}

\title{Non-coercive Neumann boundary control problems\thanks{The second author was partially supported by MCIN/ AEI/10.13039/501100011033 under research project PID2020-114837GB-I00.}}

\author{Thomas Apel\thanks{Institute of Mathematics and Computer-Based Simulation. Universit\"at der Bundeswehr M\"unchen, 85577 Neubiberg, Germany, {\tt thomas.apel@unibw.de}. \url{http://orcid.org/0000-0003-3642-3956}}
\and Mariano Mateos\thanks{Departamento de Matem\'{a}ticas, Campus de Gij\'on, Universidad de Oviedo, 33203, Gij\'on, Spain, {\tt mmateos@uniovi.es}.
\url{http://orcid.org/0000-0003-3100-412X}}
\and Arnd R\"osch\thanks{Fakult\"at f\"ur Mathematik, Universt\"at Duisburg-Essen, D-45127 Essen, Germany, {\tt arnd.roesch@uni-due.de.}
\url{http://orcid.org/0009-0001-2163-7153}}
}

\begin{document}
\maketitle

\begin{abstract}
The article examines a linear-quadratic Neumann control problem that is governed by a non-coercive elliptic equation. Due to the non-self-adjoint nature of the linear control-to-state operator, it is necessary to independently study both the state and adjoint state equations. The article establishes the existence and uniqueness of solutions for both equations, with minimal assumptions made about the problem's data. Next, the regularity of these solutions is studied in three frameworks: Hilbert-Sobolev spaces, Sobolev-Slobodecki\u\i{} spaces, and weighted Sobolev spaces. These regularity results enable a numerical analysis of the finite element approximation of both the state and adjoint state equations. The results cover both convex and non-convex domains and quasi-uniform and graded meshes. Finally, the optimal control problem is analyzed and discretized. Existence and uniqueness of the solution, first-order optimality conditions, and error estimates for the finite element approximation of the control are obtained. Numerical experiments confirming these results are included.
A significant highlight is that the discretization error estimates known from the literature, are improved even for the coercive case.
\end{abstract}

\begin{quote}
\end{quote}
\begin{quote}
\end{quote}

\section{Introduction}
\label{S1}
Let us consider a domain $\Omega\subset\mathbb R^2$ with a polygonal boundary $\Gamma$. We are concerned with the Neumann boundary control problem
\[\label{problemP}
\Pb\quad \min_{u \in \uad} J(u) := \frac{1}{2}\int_\Omega (y_u(x) - y_d(x))^2\, \dx + \frac{\nu}{2}\int_\Gamma u^2(x)\, \dx +\int_\Gamma y_u(x) g_\varphi(x)\, \dx
\]
where $y_d \in L^2(\Omega)$  and $g_\varphi\in L^2(\Gamma)$ are given functions, $\nu > 0$,
\[
\uad = \{u \in L^2(\Omega) : \umin \le u(x) \le \umax \text{ for a.{e.} } x \in \Omega\}
\]
with $-\infty \le \umin < \umax \le +\infty$, and $y_u$ is the solution of
\begin{equation}\label{E1.1}
\left\{\begin{array}{rcl} Ay + b(x)\cdot\nabla y + a_0(x) y &=& 0 \text{ in } \Omega,\\ \partial_{n_A} y &=& u\text{ on } \Gamma.\end{array}\right.
\end{equation}
Assumptions regarding the symmetric second order differential operator $A$ and the coefficients $b$ and $a_0$ will be described later. Let us just emphasize now that we will make no assumptions on $b$ and $a_0$ that would imply coerciveness of the associated bilinear form.

The main objective of this paper is to discretize the optimal control problem using the finite element method and to obtain error estimates for the approximations of the optimal control in terms of the discretization parameter $h$. The paper aims to minimize assumptions to better capture their essence. The results are valid for possibly non-convex domains and both quasi-uniform and graded meshes. Although the theory for Neumann boundary optimal control problems governed by elliptic equations is quite complete, to the best of our knowledge, the issues that arise when the elliptic operator governing the equation is not coercive have not been addressed yet; see Casas, Mateos and Tr\"oltzsch 2005 \cite{CMT05}, Casas and Mateos 2007 \cite{Casas-Mateos2007}, Mateos and R\"osch 2011 \cite{Mateos-Rosch2011}, Apel Pfefferer and R\"osch 2012 and 2015 \cite{APR2012,APR2015}, Krumbiegel and Pfefferer 2015 \cite{Krumbiegel-Pfefferer2015}. The only papers, we are aware, that  deal with optimal control problems governed by a non-coercive elliptic equation are about distributed controls; see Casas, Mateos and R\"osch 2020 and 2021 \cite{CMR2020,CMR2021}. In both papers this fact and the convexity of the domain are used in an essential way in some of the proofs, and hence those results are not applicable to our problem.

We will see that problem \Pb has a unique solution $\bar u$, and that it satisfies the optimality conditions which we state now in an informal way: there exists $\bar y$ and $\bar\varphi$ such that
\begin{subequations}
\begin{equation}\label{E1.2a}
\left\{\begin{array}{rcl} A\bar y + b(x)\cdot\nabla \bar y + a_0(x) \bar y &=& 0 \text{ in } \Omega,\\ \partial_{n_A} \bar y &=& \bar u\text{ on } \Gamma,\end{array}\right.
\end{equation}
\begin{equation}\label{E1.2b}
\left\{\begin{array}{rccl} A\bar \varphi -\nabla\cdot(b(x)\bar \varphi) + a_0(x) \bar \varphi &=& \bar y-y_d & \text{ in } \Omega,\\ \partial_{n_A} \bar \varphi + \varphi b\cdot n &=& g_\varphi & \text{ on } \Gamma,\end{array}\right.
\end{equation}
\begin{equation}\label{E1.2c}
  \int_\Gamma(\bar\varphi + \nu \bar u)(u-\bar u)\,\dx \geq 0\ \forall u\in\uad.
\end{equation}
\end{subequations}
Since \Pb is a linear-quadratic strictly convex problem, existence and uniqueness of the solution follow in a standard way once we have proved existence and uniqueness of solution of the state equation and continuity of the control-to-state mapping.
But, since we will not formulate any assumptions on $b$ or $a_0$ that would lead to a coercive operator, this task is not standard. In particular, $\mathrm{div}\,b$ may be large, such that the usual assumption $a_0-\frac12\mathrm{div}\,b\ge c_0>0$ is not satisfied. This will be done in Section \ref{S2}.

In Section \ref{S3} we investigate the regularity properties of the solutions of
the state equation and the adjoint state equations.
Since these are different, we perform this task in two steps resulting in Theorems \ref{T3.4} and \ref{T3.5}, respectively. We obtain results in Hilbert--Sobolev, in Sobolev--Slobodecki\u\i{} and in weighted Sobolev spaces, with our focus on treating the numerical approximation of \Pb in non-convex domains.
The regularity results in non-weighted spaces serve us as intermediate results to prove the error estimates in weighted Sobolev spaces, but they are also of independent interest.
Note that, although regularity results for elliptic boundary value problems are widely investigated, see, e.\,g., the monographs \cite{Grisvard85,KozlovMazyaRossmann1997,MazyaPlamenevski1984,MazyaRossmann2010,NazarovPlamenevsky1994}, the particular results which we need for our approximation theory, were not available for non-coercive problems with variable coefficients.

In Section \ref{S4} we study the numerical discretization of both the state and adjoint state equation. We obtain existence and uniqueness of the solution as well as error estimates. Our results are valid in convex and non-convex domains and for quasi-uniform and graded meshes, with possibly a non-optimal grading parameter $\mu$.

With these results at hand, we will be able to deduce existence, uniqueness, and optimality conditions  in Section \ref{S5}. Moreover, regularity properties of the optimal solution and its related state and adjoint state are given in terms of weighted Sobolev spaces. Finally, we will discretize the control problem. The control is approximated using piecewise constant functions wheras the state and adjoint state are discretized by continuous piecewise linear functions. A close inspection of the proofs in the above mentioned papers about Neumann control problems, suggests that, if no postprocessing step is done, the order of convergence of the error in $L^2(\Gamma)$ for the control variable will be limited by the order of convergence of the finite element error in $H^1(\Omega)$ for the state or the adjoint state equation; see e.g. the proof of Lemma 4.7 in \cite{CMT05}. This means that, for a non-convex domain and a quasi-uniform mesh, the order of convergence that can be obtained---applying the usual techniques in optimal control together with the regularity results and the finite element error estimates provided in this paper---is approximately $h^\lambda$, where $1/2<\lambda<1$. For instance, in the problem shown as an example in Section \ref{S6}, $h^{2/3}$ would be expected. Nevertheless, the numerical experiments show clearly order $h$, and we are able to get that in Theorem \ref{T5.7}: If the corner singularities are of type $r^{\lambda_j}$, the index $j$ counting the corners, and the mesh is graded near the corners with parameter $\mu_j$, then the approximation order of the control is $s^*\le1$ with $s^*<\frac{3\lambda_j}{2\mu_j}$, i.\,e., $s^*=1$ is achieved if $\mu_j<\frac32\lambda_j$ for all $j$, whereas Apel, Pfefferer and R\"osch \cite{APR2012,APR2015} needed the stronger grading $\mu_j<\lambda_j$. In this way, the paper does not only extend results to non-coercive problems but also improves a result for problems with coercive state equation.

\section{Existence, uniqueness and continuous dependence of the solution of the state and adjoint state equations}\label{S2}

On $A$, $b$ and $a_0$ we make the following assumptions.

\begin{assumption}\label{A2.1}
$A$ is the operator given by
\[
Ay= -  \sum_{i,k = 1}^{2}\partial_{x_k}(a_{ik}(x)\partial_{x_i}y) \ \text{ with }\ a_{ik} \in L^\infty(\Omega),
\]
$a_{ik}=a_{ki}$ for $1\leq i,k\leq 2$, and satisfying the following ellipticity condition:
\begin{equation}\label{E2.1}
\exists \Lambda > 0 \text{ such that } \sum_{i,k = 1}^{2}a_{ik}(x)\xi_i\xi_k \ge \Lambda\vert \xi \vert^2\ \ \forall \xi \in \mathbb{R}^2 \text{ and for a.a. } x \in \Omega.
\end{equation}
The function  $b:\Omega\to\mathbb{R}^2$ satisfies  $b\in L^\ps(\Omega)^2$ with $\ps > 2$ and there exists $\qs>1$ such that  $\nabla\cdot b\in L^\qs(\Omega)$ and $b\cdot n\in L^\qs(\Gamma)$. For the function $a_0:\Omega\to\mathbb{R}$ it is assumed that $a_0 \in L^\qs(\Omega)$,
$a_0(x)\geq 0$ for a.e. $x\in\Omega$ and there exists  $E\subset\Omega$ with $\vert E\vert>0$ such that $a_0(x)\geq \Lambda/2$ for all $x\in E$.
\end{assumption}

\begin{remark}\label{R2.2}
  Note that this assumption does not lead to a coercive bilinear form.
\end{remark}
Before addressing the main results of this section, we recall some well known inequalities that will be used throughout this paper.

We will often use the following form of H\"older's inequality: for $q,p_1,\cdots,p_k\in[1,\infty]$ such that $1/p_1+\ldots+1/p_k\leq 1/q$ and $f_i\in L^{p_i}(\Omega)$, $i=1,\ldots,k$ there exists a constant $C_\Omega=\vert \Omega\vert^{1/q-(1/p_1+\cdots+1/p_k)}$, such that $\Vert f_1\cdots f_k\Vert_{L^q(\Omega)}\leq C_\Omega \Vert f_1\Vert_{L^{p_1}(\Omega)}\cdots \Vert f_k\Vert_{L^{p_k}(\Omega)}$.

The inequality
\begin{equation}
\label{E2.2}\Vert y\Vert_{H^1(\Omega)} \le C_E(\Vert \nabla y\Vert_{L^2(\Omega)}+\Vert y\Vert_{L^2(E)})\quad \forall y \in H^1(\Omega)
\end{equation}
is a generalization of Poincar\'e's inequality and can be found, e.g.,  in \cite[Theorem 11.19]{Casas92B}.
In dimension 2, Sobolev's embedding theorem gives that for all $r<\infty$ there exists $K_{\Omega,r}>0$ such that
\begin{equation}
\label{E2.3} \Vert y\Vert_{L^r(\Omega)} \le K_{\Omega,r}\Vert y\Vert_{H^1(\Omega)}\quad \forall y \in H^1(\Omega).
\end{equation}

We will denote by $\langle\cdot,\cdot\rangle_\Omega$ the duality product in $H^1(\Omega)'\times H^1(\Omega)$ and by $\langle\cdot,\cdot\rangle_\Gamma$ the duality product in $H^{1/2}(\Gamma)'\times H^{1/2}(\Gamma)$. We notice that any $g\in H^{1/2}(\Gamma)'$ defines an element in $H^1(\Omega)'$, which will be denoted in the same way by
\begin{equation}\label{E2.4}\langle g,z\rangle_\Omega = \langle g,\mathrm{tr}z\rangle_\Gamma\ \forall z\in H^1(\Omega).\end{equation}
In this case, we will simply write $\langle g,z\rangle_\Gamma$. Also notice that for any fixed $q>1$, the functions $f\in L^q(\Omega)$ and $g\in L^q(\Gamma)$ define elements in $H^1(\Omega)'$ and $H^{1/2}(\Gamma)'$ respectively by
\begin{equation}\label{E2.5}\langle f,z\rangle_\Omega = \int_\Omega fz\dx,\qquad\langle g,z\rangle_\Gamma = \int_\Gamma gz\,\dx \ \forall z\in H^1(\Omega).\end{equation}

For every $y\in H^1(\Omega)$, we define $\mathcal{A}y$ by
\begin{equation}
\langle \mathcal{A}y,z\rangle_\Omega = \int_\Omega \sum_{i,k = 1}^2 a_{ik}\partial_{x_i}y\partial_{x_k} z \dx + \int_\Omega ( b\cdot\nabla y )z\dx + \int_\Omega a_0 y z\dx\ \forall z\in H^1(\Omega).
\label{E2.6}
\end{equation}
Using this operator, we have that the weak form of the state equation \eqref{E1.1} is: find $y_u\in H^1(\Omega)$ such that
\begin{equation}\label{E2.7}\langle \mathcal{A}y_u,z\rangle_\Omega = \langle u, z\rangle_\Gamma\ \forall z\in H^1(\Omega).\end{equation}
We first prove continuity of the operator $\mathcal{A}$ and G\r{a}rding's inequality. We adapt the proof of \cite[Lemma 2.1]{CMR2020}
\begin{lemma}\label{L2.3}
Under Assumption \ref{A2.1} we have that
$\mathcal{A}\in \mathcal{L}(H^1(\Omega),H^{1}(\Omega)')$ and there exists a constant $C_{\Lambda,E,b}$ such that
\begin{equation}\label{E2.8}
\langle \mathcal{A} z,z\rangle_\Omega  \ge \frac{\Lambda}{8 C_E^2}\Vert z\Vert^2_{H^1(\Omega)} - C_{\Lambda,E,b}\Vert z\Vert^2_{L^2(\Omega)}\quad \forall z \in H^1(\Omega)
\end{equation}
where $\Lambda$ and $C_E$ are the constants from \eqref{E2.1} and \eqref{E2.2}, respectively.
\end{lemma}

\begin{proof}
Let us show that $\mathcal{A}$ is a linear continuous operator. Denote $S =\{z\in H^1(\Omega):\ \Vert z\Vert_{H^1(\Omega)} = 1\}$. We split $\mathcal{A}$ into three parts $\mathcal A_i$, $i=1,2,3$.
\begin{align*}
  \Vert \mathcal A_1 y\Vert_{H^1(\Omega)'} &=   \sup_{z\in S} \int_\Omega \sum_{i,k = 1}^2 a_{ik}\partial_{x_i}y\partial_{x_k} z \dx\\
  &\leq  \sup_{z\in S} 4 \max_{1\leq i,k\leq 2}\Vert a_{ik}\Vert_{L^\infty(\Omega)} \Vert \nabla y\Vert_{L^2(\Omega)^2}\Vert \nabla z\Vert_{L^2(\Omega)^2} \\
  &\leq  4 \max_{1\leq i,k\leq 2}\Vert a_{ik}\Vert_{L^\infty(\Omega)}\Vert y\Vert_{H^1(\Omega)}.
\end{align*}

Take now $s_p>1$ such that $1/s_p = 1/\ps + 1/2$. From \eqref{E2.3} and H\"older's inequality we infer for every $y \in H^1(\Omega)$
\begin{align*}
  \Vert \mathcal A_2 y\Vert_{H^1(\Omega)'} &=   \sup_{z\in S} \int_\Omega ( b\cdot\nabla y )z\dx \leq \Vert b\cdot\nabla y\Vert_{L^{s_p}(\Omega)} \Vert z\Vert_{L^{s'_p}(\Omega)} \nonumber \\
&\le  K_{\Omega,s'_p} \Vert b\Vert_{L^\ps(\Omega)^2}\Vert \nabla y\Vert_{L^2(\Omega)^2} \leq K_{\Omega,s'_p} \Vert b\Vert_{L^\ps(\Omega)^2}\Vert y\Vert_{H^1(\Omega)},
\end{align*}
Fix now some $s_q\in(1,\qs)$ and take $r\in(1,+\infty)$ such that $1/\qs+1/r = 1/s_q$. From \eqref{E2.3} we infer that
\begin{align*}
\Vert \mathcal A_3 y\Vert_{H^1(\Omega)'} &=   \sup_{z\in S} \int_\Omega a_0 y z\dx \leq  \Vert a_0 y\Vert_{L^{s_q}(\Omega)}\Vert z\Vert_{L^{s'_q}(\Omega)} \notag \\
&\le  K_{\Omega,s'_q} \Vert a_0\Vert_{L^{\qs}(\Omega)}\Vert y\Vert_{L^{r}(\Omega)} \le K_{\Omega,s'_q} K_{\Omega,r}\Vert a_0\Vert_{L^{\qs}(\Omega)}\Vert y\Vert_{H^1(\Omega)}.
\end{align*}
Hence, we have that $\mathcal{A}$ is a well-posed linear and continuous operator.

Let us prove \eqref{E2.8}. Using  Assumption \ref{A2.1}, \eqref{E2.2}, and Young and H\"older inequalities we get
\begin{align*}
\langle \mathcal{A} z,z\rangle_\Omega &\ge \Lambda\Vert \nabla z\Vert^2_{L^2(\Omega)^2} + \frac{\Lambda}{2} \Vert z\Vert^2_{L^2(E)} - \Vert \nabla z\Vert_{L^2(\Omega)^2}\Vert bz\Vert_{L^2(\Omega)^2} \\
&\ge  \frac{\Lambda}{2}\Vert \nabla z\Vert^2_{L^2(\Omega)^2} + \frac{\Lambda}{2} \Vert z\Vert^2_{L^2(E)} - \frac{1}{2\Lambda}\Vert bz\Vert^2_{L^2(\Omega)^2}\\
 &\ge \frac{ \Lambda}{4C_E^2}\Vert z\Vert^2_{H^1(\Omega)} - \frac{1}{2\Lambda}\Vert b\Vert^2_{L^{\ps}(\Omega)^2}\Vert z\Vert^2_{L^{\frac{2 \ps}{\ps - 2}}(\Omega)}.
\end{align*}
Observe that the assumption $\ps > 2$ implies that $2 < \dfrac{2 \ps}{\ps - 2} < \infty$. Now, we apply Lions' Lemma, \cite[Chapter 2, Lemma 6.1]{Necas67}, to the chain of embeddings $H^1(\Omega) \subset\subset L^{\frac{2 \ps}{ \ps - 2}}(\Omega) \subset L^2(\Omega)$, the first one being compact and the second one continuous, to deduce the existence of a constant $C_0$ depending on $\Lambda$, $C_E$ and $\Vert b\Vert_{L^{p}(\Omega)^2}$ such that
\[
\Vert z\Vert_{L^{\frac{2\ps}{\ps - 2}}(\Omega)} \le \frac{\Lambda}{2^{3/2}\Vert b\Vert_{L^{\ps}(\Omega)^2}C_E}\Vert z\Vert_{H^1(\Omega)} + C_0\Vert z\Vert_{L^2(\Omega)}.
\]
From the last two inequalities we conclude \eqref{E2.8}
with
\[
C_{\Lambda,E,b} = \frac{C_0^2\Vert b\Vert^2_{L^{\ps}(\Omega)^2}}{\Lambda}
\]
and the proof is complete.
\end{proof}

\begin{remark}\label{R2.4}
  Notice that, to prove Lemma \ref{L2.3}, we use neither $\nabla\cdot b\in L^\qs(\Omega)$ nor $b\cdot n\in L^\qs(\Gamma)$ for some $\qs>1$.
\end{remark}

The adjoint operator of $\mathcal A$ is $\mathcal{A}^*$. We have $\mathcal{A}^*z\in H^1(\Omega)'$ for every $z\in H^1(\Omega)$. In the next lemma, we justify that under the mild Assumption \ref{A2.1}, we can integrate by parts and use the expected form of the adjoint equation \eqref{E1.2b}.
\begin{lemma}\label{L2.5} Suppose that Assumption \ref{A2.1} holds. Then
\[\langle \mathcal{A}^* z ,y \rangle_\Omega = \int_\Omega \sum_{i,k = 1}^2 a_{ki}\partial_{x_i}z\partial_{x_k} y \dx - \int_\Omega y\nabla\cdot ( b z)\dx +
\int_\Gamma y z b\cdot n\,\dx
+ \int_\Omega a_0 y z\dx.
\]
for all $y\in H^1(\Omega)$.
\end{lemma}
\begin{proof}
By definition
\[\langle \mathcal{A}^* z ,y \rangle_\Omega=\langle \mathcal{A} y ,z \rangle_\Omega \ \forall y,z\in H^1(\Omega)\]
and we only have to justify that, under Assumption \ref{A2.1}, we can do integration by parts to get
\[\int_\Omega (b\cdot\nabla y) z\dx = -\int_\Omega y\nabla \cdot(bz)\dx + \int_\Gamma yz b\cdot n \,\dx.\]
This is equivalent to proving that we can apply the Gauss theorem to obtain
\[\int_\Omega \nabla\cdot (yzb)\dx = \int_\Gamma yz b\cdot n \,\dx.\]
Using that $y,z\in H^1(\Omega)\hookrightarrow L^r(\Omega)$ for all $r<+\infty$, $b\in L^\ps(\Omega)^2$ for some $\ps>2$ and $\nabla\cdot b\in L^\qs(\Omega)$ for some $\qs>1$, applying H\"older's inequality, we have
\[\nabla(yz)\cdot b=z\nabla y\cdot b+ y\nabla z\cdot b\in L^{\frac{2\ps}{2+\ps}}(\Omega)\text{ and }yz\nabla\cdot b\in L^{\frac{\qs+1}{2}}(\Omega),\]
so $\nabla\cdot(yzb)\in L^s(\Omega)$, where $s = \min\left\{ \dfrac{2\ps}{2+\ps},\dfrac{\qs+1}{2}\right\}$ satisfies $1 <s <2$. From Assumption \ref{A2.1}, it is also clear that $yzb\in L^s(\Omega)^2$, and using \cite[Lema II.1.2.2]{Sohr2001}, we deduce that  $yzb$ has a normal trace $yzb\cdot n$ defined in the space of $(W^{1-1/s',s'}(\Gamma))'$ {\em via} Gauss theorem: for every $v\in W^{1,s'}(\Omega)$
\[\langle yzb\cdot n,v\rangle_{(W^{1-1/s',s'}(\Gamma))',W^{1-1/s',s'}(\Gamma)} = \int_\Omega \nabla\cdot (v yzb) \dx.\]
Since we are assuming that $b\cdot n\in L^\qs(\Gamma)$ for some $\qs>1$, then $yz b\cdot n\in L^{\frac{\qs+1}{2}}(\Gamma)\hookrightarrow L^s(\Gamma)$. Therefore, we have that
\[\langle yzb\cdot n,v\rangle_{(W^{1-1/s',s'}(\Gamma))',W^{1-1/s',s'}(\Gamma)} = \int_\Gamma v yzb\cdot n\,\dx.\]
Taking $v = 1$ in the above equalities, we have that
\[\int_\Omega \nabla\cdot (yzb) \dx = \int_\Gamma yzb\cdot n\,\dx,\]
and the proof is complete.
\end{proof}
Next, we adapt the proof of \cite[Theorem 2.2]{CMR2020} to show existence and uniqueness of the solution of the state equation.
\begin{lemma}\label{L2.6}
Under Assumption \ref{A2.1}, the linear operator $\mathcal{A}:H^1(\Omega) \longrightarrow H^{1}(\Omega)'$ is an isomorphism.
\end{lemma}
\begin{proof}
  Let us first see that $\mathcal{A}$ is injective. Consider $y\in H^1(\Omega)$ such that $\mathcal{A}y =0$. We will prove that $y\leq 0$. The contrary inequality follows by arguing on $-y$. Suppose there exists some $\mathcal{O}\subset\Omega$ with positive measure such that $y (x) > 0$ if $x\in \mathcal{O}$. Take $0 < \rho < \text{ess}\sup_{x\in\Omega} y(x)\leq +\infty$ and define $y_\rho(x) = (y(x)-\rho)^+ = \max \{y(x)-\rho,0\}$. Denote $\Omega_\rho = \{x\in\Omega:\ \nabla y_\rho(x)\neq 0\}$.
   Notice that $y_\rho\in H^1(\Omega)$,
  \[\nabla y_\rho(x) = \left\{\begin{array}{cc}
                                \nabla y(x) & \text{ if }y(x) > \rho \\
                                0  & \text{ if }y(x) \leq \rho,
                              \end{array}
  \right.
  \]
  which means that $\Omega_\rho\subset\{x:\ y(x)>\rho\}$. We also remark that $y_\rho(x) = 0\text{ if }y(x)\leq 0$, and that $y(x)\geq y_\rho(x) \geq  0$ if $y(x) \geq 0 $.
  Using these properties, and H\"older's and Young's inequalities, we have that
  \begin{align*}
    0 &=  \langle \mathcal{A} y ,y_\rho \rangle_\Omega = \int_\Omega \sum_{i,k = 1}^2 a_{ik}\partial_{x_i}y\partial_{x_k} y_\rho  \dx + \int_\Omega ( b\cdot\nabla y )y_\rho \dx + \int_\Omega a_0 y y_\rho \dx\\
    &\geq  \int_{\Omega_\rho} \sum_{i,k = 1}^2 a_{ik}\partial_{x_i}y_\rho\partial_{x_k} y_\rho  \dx + \int_{\Omega_\rho} ( b\cdot\nabla y_\rho )y_\rho \dx + \int_{\Omega} a_0 y_\rho y_\rho \dx\\
    &\geq  \Lambda \Vert \nabla y_\rho\Vert^2_{L^2(\Omega_\rho)} - \Vert b\Vert_{L^\ps(\Omega_\rho)^2}\Vert \nabla y_\rho\Vert_{L^2(\Omega_\rho)} \Vert y_\rho\Vert_{L^\frac{2\ps}{\ps-2}(\Omega_\rho)} + \frac{\Lambda}{2}\Vert y_\rho\Vert_{L^2(E)}^2\\
    &\geq \frac{\Lambda}{2}\Vert \nabla y_\rho\Vert^2_{L^2(\Omega_\rho)} -\frac{1}{2\Lambda} \Vert b\Vert_{L^\ps(\Omega_\rho)^2}^2 \Vert y_\rho\Vert_{L^\frac{2\ps}{\ps-2}(\Omega_\rho)}^2+
    \frac{\Lambda}{2}\Vert y_\rho\Vert_{L^2(E)}^2\\
    & =   \frac{\Lambda}{2}\Vert \nabla y_\rho\Vert^2_{L^2(\Omega)} -\frac{1}{2\Lambda} \Vert b\Vert_{L^\ps(\Omega_\rho)^2}^2 \Vert y_\rho\Vert_{L^\frac{2\ps}{\ps-2}(\Omega_\rho)}^2+
    \frac{\Lambda}{2}\Vert y_\rho\Vert_{L^2(E)}^2
  \end{align*}
  Next we use that $\Omega_\rho\subset \Omega$, \eqref{E2.3}, \eqref{E2.2} and the just proved inequality to obtain:
  \begin{align*}
    \Vert y_\rho\Vert_{L^\frac{2\ps}{\ps-2}(\Omega_\rho)}^2  &\leq  \Vert y_\rho\Vert_{L^\frac{2\ps}{\ps-2}(\Omega)}^2 \leq K_{\Omega,\frac{2p}{p-2}}^2 \Vert y_\rho\Vert_{H^1(\Omega)}^2 \\
    &\leq  2 K_{\Omega,\frac{2\ps}{\ps-2}}^2 C_E^2 \left(\Vert \nabla y_\rho\Vert^2_{L^2(\Omega)} + \Vert y_\rho\Vert_{L^2(E)}^2\right)\\
    &\leq \frac{2 K_{\Omega,\frac{2\ps}{\ps-2}}^2 C_E^2}{\Lambda^2} \Vert b\Vert_{L^\ps(\Omega_\rho)^2}^2 \Vert y_\rho\Vert_{L^\frac{2\ps}{\ps-2}(\Omega_\rho)}^2
  \end{align*}
  We can deduce from this a positive lower bound for the norm of $b$ in $L^\ps(\Omega_\rho)^2$ independent of $\rho$.
  \[\Vert b\Vert_{L^\ps(\Omega_\rho)^2}\geq \frac{\Lambda}{\sqrt{2}K_{\Omega,\frac{2\ps}{\ps-2}} C_E} > 0.\]
  But we have that $\vert \Omega_\rho\vert\to 0$ as $\rho\to  \text{ess}\sup_{x\in\Omega} y(x)$; see \cite[Theorem 2.2]{CMR2020}. So we have achieved a contradiction.

  \smallskip

  Finally we have just to check that the range of $\mathcal{A}$ is dense and closed. Since we already have established G\r{a}rding's inequality \eqref{E2.8} for the operator $\mathcal{A}$, the proof of closeness done in \cite[Theorem 2.2]{CMR2020} applies to our case changing the norms in $H^1_0(\Omega)$ and $H^{-1}(\Omega)$ respectively by the norms in $H^1(\Omega)$ and its dual space, and thus it is omitted. By a well known duality argument, the denseness of the range of  $\mathcal{A}$ follows from the injectivity of $\mathcal{A}^*$.

  The argument used above to obtain the injectivity of $\mathcal{A}$ does not work for $\mathcal{A}^*$.
  Notice that at one moment we use that $\int_\Omega (b\cdot\nabla y) y_\rho\dx = \int_{\Omega_\rho} (b\cdot\nabla y_\rho) y_\rho\dx$. When dealing with the adjoint operator, we would find the term $\int_\Omega (b\cdot \nabla z_\rho) z\dx$, which in general is different from $\int_{\Omega_\rho} (b\cdot \nabla z_\rho) z_\rho\dx$. But we can
  obtain injectivity of the adjoint operator as follows.
  Consider $z\in H^1(\Omega)$ such that $\mathcal{A}^*z =0$. For all $\varepsilon\geq 0$ define
  \[\Omega^\varepsilon = \{x\in\Omega:\ \vert z(x)\vert > \varepsilon\}\]
  Let us see that $\vert \Omega^0 \vert = 0$, which readily implies that $z=0$. Let us define $z^\varepsilon(x) = \proj_{[-\varepsilon,\varepsilon]}(z(x))$. Using integration by parts, that $z=0$ in $\Omega\setminus\Omega^0$, that $\nabla z^\varepsilon =0$ in $\Omega^\varepsilon$ and $\nabla z^\varepsilon = \nabla z$ in $\Omega\setminus\Omega^\varepsilon$, and that $z z^\varepsilon \geq (z^\varepsilon)^2$, we have
  \begin{align*}
    0 &=  \langle \mathcal{A}^*z ,z^\varepsilon \rangle_\Omega \\
    &=   \int_\Omega \sum_{i,k = 1}^2 a_{ki}\partial_{x_i}z\partial_{x_k} z^\varepsilon \dx - \int_\Omega z^\varepsilon\nabla\cdot ( b z)\dx +
\int_\Gamma z^\varepsilon z b\cdot n\,\dx
+ \int_\Omega a_0 z^\varepsilon z\dx\\
 &=  \int_\Omega \sum_{i,k = 1}^2 a_{ki}\partial_{x_i}z\partial_{x_k} z^\varepsilon \dx + \int_\Omega z b\cdot\nabla z^\varepsilon\dx
+ \int_\Omega a_0 z^\varepsilon z\dx\\
 &\geq  \Lambda \Vert \nabla z^\varepsilon\Vert_{L^2(\Omega)^2} - \Vert b\Vert_{L^{\ps}(\Omega^0\setminus\Omega^\varepsilon)} \Vert \nabla z^\varepsilon\Vert_{L^2(\Omega)^2} \Vert z^\varepsilon\Vert_{L^{\frac{2\ps}{\ps-2}}(\Omega^0\setminus\Omega^\varepsilon)} +\frac{\Lambda}{2}\Vert z^\varepsilon\Vert_{L^2(E)}^2\\
 &\geq  \frac{\Lambda}{2} \Vert \nabla z^\varepsilon\Vert_{L^2(\Omega)^2} -\frac{1}{2\Lambda}
 \Vert b\Vert_{L^{\ps}(\Omega^0\setminus\Omega^\varepsilon)}^2 \Vert z^\varepsilon\Vert_{L^{\frac{2\ps}{\ps-2}}(\Omega^0\setminus\Omega^\varepsilon)}^2
 +\frac{\Lambda}{2}\Vert z^\varepsilon\Vert_{L^2(E)}^2 \\
    \end{align*}

 So, using this and that $\vert z^\varepsilon(x)\vert\leq \varepsilon$ for a.e. $x\not\in\Omega^\varepsilon$ we get
 \begin{align*}
   \Vert z^\varepsilon\Vert_{H^1(\Omega)}^2 &\leq  2 C_E^2 \left(\Vert \nabla z^\varepsilon\Vert_{L^2(\Omega)^2} + \Vert z^\varepsilon\Vert_{L^2(E)}^2\right) \\
   &\leq \frac{2 C_E^2}{\Lambda^2} \Vert b\Vert_{L^{\ps}(\Omega^0\setminus\Omega^\varepsilon)}^2 \Vert z^\varepsilon\Vert_{L^{\frac{2\ps}{\ps-2}}(\Omega^0\setminus\Omega^\varepsilon)}^2
    \leq
   \frac{2 C_E^2}{\Lambda^2} \Vert b\Vert_{L^{\ps}(\Omega^0\setminus\Omega^\varepsilon)}^2 \vert \Omega^0\setminus\Omega^\varepsilon\vert^{\frac{\ps-2}{\ps}} \varepsilon^2.
 \end{align*}
 On the other hand, using that $\vert z^\varepsilon\vert =\varepsilon$ in $\Omega^\varepsilon$ and the previous inequality, we have
 \begin{align*}
 \vert \Omega^\varepsilon\vert  &=  \frac{1}{\varepsilon^2}\int_{\Omega^\varepsilon}z^\varepsilon(x)^2\dx \leq \frac{1}{\varepsilon^2}\Vert z^\varepsilon\Vert_{L^2(\Omega)}^2
 \leq \frac{1}{\varepsilon^2}\Vert z^\varepsilon\Vert_{H^1(\Omega)}^2 \\
 &\leq  \frac{2 C_E^2}{\Lambda^2} \Vert b\Vert_{L^{\ps}(\Omega^0\setminus\Omega^\varepsilon)}^2 \vert \Omega^0\setminus\Omega^\varepsilon\vert^{\frac{\ps-2}{\ps}}
 \end{align*}
 Since $\vert \Omega^0\setminus\Omega^\varepsilon\vert  = \mathrm{meas}\{x\in\Omega: 0 < \vert z(x)\vert < \varepsilon\}\to 0$ as $\varepsilon\to 0$, we have proved that $\vert \Omega^\varepsilon\vert\to 0$ as $\varepsilon\to 0$ and hence $\vert \Omega^0\vert =0$.
\end{proof}
\begin{corollary}\label{C2.7}
Under Assumption \ref{A2.1}, the linear operator $\mathcal{A}^*:H^1(\Omega) \longrightarrow H^{1}(\Omega)'$ is an isomorphism.
\end{corollary}

\section{Regularity of the solution of the state and adjoint state equations}\label{S3}

To obtain further regularity, from now on we will suppose

\medskip
\begin{assumption}\label{A3.1}
The coefficients $a_{ik}$ belong to $C^{0,1}(\bar\Omega)$, $1\leq i,k,\leq 2$.
\end{assumption}

Let us denote by $m$ the number of sides of $\Gamma$ and $\{\xS_j\}_{j=1}^m$ its vertices, ordered counterclockwise. For convenience denote also $\xS_0=\xS_m$ and $\xS_{m+1}=\xS_1$. We denote by $\Gamma_j$ the side of $\Gamma$ connecting $\xS_{j}$ and $\xS_{j+1}$, and by $\omega_j\in (0,2\pi)$ the angle interior to $\Omega$ at $\xS_j$, i.e., the angle defined by $\Gamma_{j}$ and $\Gamma_{j-1}$, measured counterclockwise. Notice that  $\Gamma_{0}=\Gamma_m$. We use $(r_j,\theta_j)$ as local polar coordinates at $\xS_j$, with $r_j=\vert x-\xS_j\vert$ and $\theta_j$ the angle defined by $\Gamma_j$ and the segment $[\xS_j,x]$.
In order to describe and analyze the regularity of the functions near the corners, we will introduce for every $j\in\{1,\ldots,m\}$ the infinite cone
\[K_j=\{x\in\mathbb{R}^2:0<r_j,\,0<\theta_j<\omega_j\}.\]

For every $j\in\{1,\ldots,m\}$ we call $A_j$ the operator with constant coefficients, corresponding to the corner $\xS_j$, given by
\[A_j y = \sum_{i,k = 1}^{2}\partial_{x_k}(a_{ik}(\xS_j)\partial_{x_i}y).\]
We denote by $\lambda_j$ the leading singular exponent associated with the operator $A_j$ at the corner $\xS_j$, i.e., the smallest $\lambda_j>0$ such that there exists a solution  of the form $y_j=r_j^{\lambda_j}\varphi_j(\theta_j)$, with $\varphi_j$ smooth enough, of
\[A_j y_j =0\mbox{ in }K_j,\ \partial_{n_{A_j}} y_j =0\mbox{ on }\partial K_j.\]
For instance, for $Ay = -\Delta y$ we have $\lambda_j = \pi/\omega_j$. We denote $\lambda = \min\{\lambda_j\}$.

With the usual technique of taking a partition of the unity to localize the problem in the corners, freezing the coefficients and doing an appropriate linear change of variable, the classical results for the Laplace operator are also valid in our case; see, e.g. \cite[Section 2.1]{Mateos2000} for a detailed example of application of this technique. Notice that the symmetry hypothesis $a_{ik}=a_{ki}$ introduced in Assumption \ref{A2.1} implies that the same change of variable that transforms $A_j$ into $-\Delta$ will transform the conormal derivative $\partial_{n_{A_j}}$ into the normal derivative $\partial_n$ in the new variables, and not in an oblique derivative.

We continue with regularity results for problems with $b\equiv0$ and $a_0\equiv0$ and use the standard Sobolev and Sobolev--Slobodetski\u\i{} spaces but also weighted Sobolev spaces as follows. Let $k\in\mathbb{N}_0$, $1\le p\le\infty$, and $\vec\beta=(\beta_1,\ldots,\beta_m)^T\in\mathbb{R}^m$, $j\in\{1,\ldots,m\}$. For ball-neighborhoods $\Omega_{R_j}$ of $\xS_j$ with radius $R_j\leq 1$ and $\Omega^0:=\Omega\setminus\bigcup_{j=1}^m \Omega_{R_j/2}$ we define norms via
\begin{align*}
  \Vert v\Vert_{W^{k,p}_{\beta_j}(\Omega_{R_j})}^p & = \sum_{\mid\alpha\mid\le k}
  \Vert r_j^{\beta_j}D^\alpha v\Vert_{L^p(\Omega_{R_j})}^p, \\
  \Vert v\Vert_{V^{k,p}_{\beta_j}(\Omega_{R_j})}^p & = \sum_{\mid\alpha\mid\le k}
  \Vert r_j^{\beta_j-k+\mid\alpha\mid}D^\alpha v\Vert_{L^p(\Omega_{R_j})}^p,
\end{align*}
 where the standard modification for $p=\infty$ is used. The spaces $W^{k,p}_{\vec\beta}(\Omega)$ and $V^{k,p}_{\vec\beta}(\Omega)$ denote the set of all functions $v$ such that
\begin{align*}
  \Vert v\Vert_{W^{k,p}_{\vec\beta}(\Omega)} & := \Vert v\Vert_{W^{k,p}(\Omega^0)} + \sum_{j=1}^m
  \Vert v\Vert_{W^{k,p}_{\beta_j}(\Omega_{R_j})}, \\
  \Vert v\Vert_{V^{k,p}_{\vec\beta}(\Omega)} & := \Vert v\Vert_{W^{k,p}(\Omega^0)} + \sum_{j=1}^m
  \Vert v\Vert_{V^{k,p}_{\beta_j}(\Omega_{R_j})},
\end{align*}
respectively, are finite. The corresponding seminorms are defined by setting $\vert \alpha\vert=k$ instead of $\vert\alpha\vert\le k$. For the definition of the corresponding trace spaces $W^{k-1/p,p}_{\vec\beta}(\Gamma_j)$, $V^{k-1/p,p}_{\vec\beta}(\Gamma_j)$, $W^{k-1/p,p}_{\vec\beta}(\Gamma)$ and $V^{k-1/p,p}_{\vec\beta}(\Gamma)$ we refer to \cite[Sect.~6.2.10]{MazyaRossmann2010}, see also \cite[Section 2.2]{Pfefferer2014}.
We will also use the notation $L^p_{\vec\beta}(\Omega)$ for $W^{0,p}_{\vec\beta}(\Omega)$.

\begin{lemma}\label{L3.2}
Suppose that Assumption \ref{A3.1} holds. Consider $f\in H^1(\Omega)'$ and $g\in H^{1/2}(\Gamma)'$ such that
\[\langle f,1\rangle_\Omega + \langle g,1\rangle_\Gamma = 0,\]
and let $y\in H^{1}(\Omega)$ be the unique solution, up to a constant, of
\[\int_\Omega \sum_{i,k=1}^2 a_{ik}\partial_{x_i}y\partial_{x_k} z\dx = \langle f,z\rangle_{\Omega} + \langle g,z\rangle_\Gamma \quad \forall z\in H^1(\Omega).\]
We have the following regularity results.

\hspace{-\parindent}\textup{(a)} If $f\in H^{2-t}(\Omega)'$, and $g\in\displaystyle\prod_{j=1}^m H^{t-3/2}(\Gamma_j)$ for some $1<t<1+\lambda$,
   $t\leq 2$, then
\[y\in H^{t}(\Omega)\mbox{ and }\vert y\vert_{H^{t}(\Omega)}\leq C_{A,t}\Big(\Vert f\Vert_{ H^{2-t}(\Omega)'}+\sum_{j=1}^m\Vert g\Vert_{H^{t-3/2}(\Gamma_j)}\Big).\]

\hspace{-\parindent}\textup{(b)} If $f\in L^r(\Omega)$ and $g\in \displaystyle\prod_{j=1}^m W^{1-1/r,r}(\Gamma_j)$ for some $1<r<\dfrac{2}{2-\lambda}$ if $\lambda<2$, $r>1$ arbitrary if $\lambda\ge2$, then
\[y\in W^{2,r}(\Omega)\mbox{ and }\vert y\vert_{W^{2,r}(\Omega)}\leq C_{A,r}\Big(\Vert f\Vert_{ L^{r}(\Omega)}+\sum_{j=1}^m\Vert g\Vert_{W^{1-1/r,r}(\Gamma_j)}\Big).\]

\hspace{-\parindent}\textup{(c)} Consider $s\in(1,\infty)$ and ${\vec\beta}$ such that  $2-\dfrac{2}{s}-\lambda_j<\beta_j<2-\dfrac2s$, $\beta_j\geq 0$ for all $j\in\{1,\ldots,m\}$. If $f\in L^s_{\vec\beta}(\Omega)$ and $g\in \prod_{j=1}^m V^{1-1/s,s}_{\vec\beta}(\Gamma_j)$,  then
\[y\in W^{2,s}_{\vec\beta}(\Omega)\mbox{ and }\vert y\vert_{W^{2,s}_{\vec\beta}(\Omega)}\leq C_{A,{\vec\beta}}\Big(\Vert f\Vert_{ L^{s}_{\vec\beta}(\Omega)}+\sum_{j=1}^m\Vert g\Vert_{V^{1-1/s,s}_{\vec\beta}(\Gamma_j)}\Big).\]
\end{lemma}
\begin{remark}\label{R3.3}

Let us briefly comment on the function spaces appearing in the lemma.
Notice that for $t=2$, $H^{t-2}(\Omega) = H^{2-t}(\Omega)=L^2(\Omega)$, and for $3/2<t$, $H^{2-t}(\Omega) = H^{2-t}_0(\Omega)$ and hence $H^{t-2}(\Omega) = H^{2-t}(\Omega)'$. Nevertheless, for $1<t<3/2$, $H^{2-t}(\Omega)' \neq H^{t-2}(\Omega)$.
Also take into account that
\begin{align*}
 \displaystyle\prod_{j=1}^m H^{t-3/2}(\Gamma_j) = H^{t-3/2}(\Gamma) & \text{ if } t<2,&
  \displaystyle\prod_{j=1}^m W^{1-1/r,r}(\Gamma_j)= W^{1-1/r,r}(\Gamma) & \text{ if } r<2
\end{align*}
We remark that the mapping  $y\mapsto\partial_{n_A}y$ is linear and continuous from $H^2(\Omega)$ onto $\prod_{j=1}^m H^{1/2}(\Gamma_j)$; see \cite[Theorem 1.5.2.8]{Grisvard85}.

Regarding weighted spaces, we notice that $V^{1-1/s,s}_{\vec\beta}(\Gamma)=W^{1-1/s,s}_{\vec\beta}(\Gamma)$ if $\beta_j>1-\frac2s$ or $\beta_j<-\frac2s$ for all $j\in\{1,\ldots,m\}$, while these spaces differ by a constant in the vicinity of each corner $\xS_j$ where $-\frac2s<\beta_j<1-\frac2s$,
see \cite[Theorem 2.1]{MazyaPlamenevski1984} or \cite[page 131]{NazarovPlamenevsky1994}.
\end{remark}
\begin{proof}[Proof of Lemma \ref{L3.2}]
The result in (a) can be deduced  from \cite[Theorem 9.2]{FabesMendezMitrea1998} for $1<t<3/2$, from \cite[Theorem (23.3)]{Dauge-1988} for $3/2<t<2$, and from \cite[Corollary 4.4.4.14]{Grisvard85} for $t=2$. The case $t=3/2$ follows by interpolation.
Statement (b) follows from \cite[Corollary 4.4.4.14]{Grisvard85}. Part (c) follows by standard arguments but we did not find this particular result in the literature. Therefore we sketch the proof here for the case of constant coefficients. As said above, the result in the case of Lipschitz coefficients follows from this one using the  localization-and-freezing technique.

We will use \cite[Theorem 1.2.5]{MazyaRossmann2010} stating a similar result for a cone $K$ and weighted $V$-spaces. For the problem under consideration and in our notation it says that $y\in V^{2,s}_\beta(K)$ if $f\in L^s_\beta(K)$ and $g\in V^{1-1/s,s}_\beta(\partial K\setminus O)$ provided that $s\in(1,\infty)$ and $2-\frac2s-\beta\not\in\{k\lambda, k\in\mathbb{Z}\}$. To satisfy the latter condition we assume $2-\frac{2}{s}-\lambda_j<\beta_j < 2-\frac2s$ for all $j\in\{1,\ldots,m\}$.

The reformulation from the vicinity of a vertex of the domain $\Omega$ is achieved by using cut-off functions $\zeta_j:\Omega\to[0,1]$ with $\zeta_j\equiv 1$ in $\Omega_{R_j/2}$, $\zeta_j\equiv 0$ in $\Omega\setminus\Omega_{R_j}$, and $\partial_{n_{A_j}}\zeta_j=0$ on $\partial\Omega\cap\partial\Omega_{R_j}$. We split $y\in H^1(\Omega)$ into
\[
  y=\sum_{j=1}^m y_j+w,\quad\text{where}\quad y_j=\zeta_j(y-y(\xS_j)).
\]
With this construction we get $y_j(\xS_j)=0$ and $\textrm{supp}\,y_j=\bar\Omega_{R_j}$ such that we can consider the problem $A y_j=f_j$ with Neumann boundary condition $\partial_{n_{A_j}}y_j=g_j=\zeta_j g$ in the cone $K_j$. For $f_j$, we have
\begin{align*}
  f_j=A\big(\zeta_j(y-y(\xS_j))\big)=\begin{cases}
    \zeta_jf & \text{in }\Omega_{R_j/2} \\
    \zeta_jf-\mathfrak{b}_j\cdot\nabla y-\mathfrak{a}_j(y-y(\xS_j)) & \text{in }\Omega_{R_j}\setminus\Omega_{R_j/2} \\
    0 & \text{in }K_j\setminus\Omega_{R_j}
  \end{cases}
\end{align*}
with smooth functions $\mathfrak{b}_j$ and $\mathfrak{a}_j$ due to the constant coefficients in $A$.
From  $f\in L^s_{\vec\beta}(\Omega)$ and $y\in H^1(\Omega)$ we conclude $f_j\in L^{\hat s}_{\beta_j}(K)$, $\hat s=\min(s,2)$ where we use that $\beta_j\ge0$. Moreover, the assumption $g\in \prod_{j=1}^mV^{1-1/s,s}_{\vec\beta}(\Gamma_j)$ leads to $g_j\in V^{1-1/s,s}_{\beta_j}(\partial K_j\setminus O_j)$ such that \cite[Theorem 1.2.5]{MazyaRossmann2010} leads to $y_j\in V^{2,\hat s}_{\beta_j}(K_j)\hookrightarrow W^{2,\hat s}_{\beta_j}(K_j)$. Since the function $w$ does not contain corner singularities, hence $w\in W^{2,s}(\Omega)$,  we obtain $y\in  W^{2,\hat s}_{\vec\beta}(\Omega)$. If $s\le2$ we are done.

Otherwise, when $s>2$, we have $y\in H^2(\Omega_{R_j}\setminus\Omega_{R_j/2})\hookrightarrow W^{1,s}(\Omega_{R_j}\setminus\Omega_{R_j/2})$, and we reiterate $f_j\in L^s_{\beta_j}(K)$ and $y_j\in V^{2,s}_{\beta_j}(K_j)\hookrightarrow W^{2,\hat s}_{\beta_j}(K_j)$ leading to $y\in  W^{2,s}_{\vec\beta}(\Omega)$.
\end{proof}

\begin{theorem}\label{T3.4}
Suppose that Assumptions \ref{A2.1} and \ref{A3.1} hold.
Consider $f\in H^1(\Omega)'$ and $u\in H^{1/2}(\Gamma)'$
and let $y\in H^{1}(\Omega)$ be the unique solution of
\begin{equation}\label{E3.1}\langle \mathcal{A}y,z\rangle_\Omega = \langle f,z\rangle_\Omega + \langle  u,z\rangle_\Gamma\quad \forall z\in H^1(\Omega).\end{equation}
We have the following regularity results.

\hspace{-\parindent}\textup{(a)} If $a_0\in L^q(\Omega)$,
$f\in H^{2-t}(\Omega)'$ and $u\in \prod_{j=1}^m H^{t-3/2}(\Gamma_j)$ for some $t$ such that $1<t<1+\lambda$, $t\leq 2$ and
$q = \dfrac{2}{3-t}$,
then $y\in H^{t}(\Omega)$  and there exists a constant $C_{\mathcal{A},t}>0$ such that
\[\Vert y\Vert_{H^{t}(\Omega)}\leq C_{\mathcal{A},t}(\Vert f\Vert_{H^{2-t}(\Omega)'}+\sum_{j=1}^m\Vert u\Vert_{H^{t-3/2}(\Gamma_j)}).\]

\hspace{-\parindent}\textup{(b)} If $a_0\in L^r(\Omega)$, $f\in L^r(\Omega)$ and $u\in\prod_{j=1}^m W^{1-1/r,r}(\Gamma_j)$ for some $r\in(1,\ps]$ satisfying $r<\dfrac{2}{2-\lambda}$ in case of $\lambda<2$, then $y\in W^{2,r}(\Omega)$ and there exists a constant $C_{\mathcal{A},r}>0$ such that
\[\Vert y\Vert_{W^{2,r}(\Omega)}\leq C_{\mathcal{A},r}(\Vert f\Vert_{ L^{r}(\Omega)}+\sum_{j=1}^m\Vert u\Vert_{ W^{1-1/r,r}(\Gamma_j)}).\]

\hspace{-\parindent}\textup{(c)} If
  $a_0\in L^{ p}_{\vec\beta}(\Omega)$, $f\in L^{p}_{\vec\beta}(\Omega)$ and $u\in \prod_{j=1}^m W^{1-1/{ p}, p}_{\vec\beta}(\Gamma_j)$ for some $p\in(1,2]$ and some $\vec\beta$ such that $2-\dfrac{2}{p}-\lambda_j<\beta_j <2-\dfrac{2}{p}$  and $\beta_j\geq 0$ for all $j\in\{1,\ldots,m\}$, then $y\in W^{2,p}_{\vec\beta}(\Omega)$ and there exists a constant $C_{\mathcal{A},\vec\beta,p}>0$ such that
\[\Vert y\Vert_{W^{2,p}_{\vec\beta}(\Omega)}\leq C_{\mathcal{A},{\vec\beta},p}(\Vert f\Vert_{ L^{p}_{\vec\beta}(\Omega)}+\sum_{j=1}^m\Vert u\Vert_{W^{1-1/{p},p}_{\vec\beta}(\Gamma_j)}).\]

\end{theorem}

\begin{proof}
 Let us define \[F = -b\cdot\nabla y -a_0 y.\] From the proof of Lemma \ref{L2.3}, we know that $F\in H^1(\Omega)'$. Also, taking $z=1$ in \eqref{E3.1}, we have that
\[\langle f+F,1\rangle_\Omega + \langle u,1\rangle_\Gamma = 0,\]
so the conditions of Lemma \ref{L3.2} apply to the problem
\[\langle A y,z\rangle_\Omega = \langle f+F,z\rangle_\Omega + \langle  u,z\rangle_\Gamma\quad \forall z\in H^1(\Omega).\]
We have to investigate the regularity of $F$.

\medskip

\hspace{-\parindent}(a) For $1<\tau\leq t$, define $S = \{z\in H^{2-\tau}(\Omega):\ \Vert z\Vert_{H^{2-\tau}(\Omega)} = 1\}$. We have that $F \in H^{2-\tau}(\Omega)'$ if and only if
\[\Vert F\Vert_{H^{2-\tau}(\Omega)'} =  \sup_{z\in S}\vert\langle F,z\rangle_\Omega\vert  < +\infty.\]
Applying H\"older's inequality, we can deduce the existence of a constant $C_\Omega>0$, that may depend on the measure of $\Omega$,  such that
\begin{align}
  \vert\langle F,z\rangle_\Omega\vert &= \left\vert\int_\Omega (b\cdot\nabla y+a_0 y)z\dx \right\vert\notag \\
  &\leq  C_\Omega\big(\Vert b\Vert_{L^\ps(\Omega)^2} \Vert \nabla y\Vert_{L^{r_p}(\Omega)} + \Vert a_0\Vert_{L^q(\Omega)}  \Vert y\Vert_{L^{r_q}(\Omega)}\big)\Vert z\Vert_{L^s(\Omega)} \label{E3.2}
\end{align}
where
\begin{equation}
\label{E3.3}
\frac{1}{\ps}+\frac{1}{r_p}+\frac{1}{s}\leq 1,\quad \frac{1}{q}+\frac{1}{r_q}+\frac{1}{s}\leq 1.
\end{equation}
Let us also notice that $H^{2-\tau}(\Omega)\hookrightarrow L^s(\Omega)$ if and only if
\begin{equation}\label{E3.4}
  \tau = 1 + \frac{2}{s}.
\end{equation}
We will apply a boot-strap argument.

\medskip

\hspace{-\parindent}Step 1. We know that $y\in H^1(\Omega)$, so $r_p=2$ and for $r_q$ we can take any real number. Noting that $q>1$, using \eqref{E3.3} and taking
\[\frac{1}{s} = \min\left\{1-\frac{1}{\ps}-\frac{1}{r_p},1-\frac{1}{q}-\frac{1}{r_q}\right\},\]
we have that $1/s>0$ for $r_q$ big enough and both conditions in \eqref{E3.3} are satisfied. Hence
we deduce that $F\in H^{2-\tau}(\Omega)'$. Since $u\in \prod_{j = 1}^{m}H^{t-3/2}(\Gamma_j)$, a direct application of Lemma \ref{L3.2} yields that $y\in H^{\min\{t,\tau\}}(\Omega)$.
If $\tau\geq t$, the proof is complete.

\medskip
\hspace{-\parindent}Step 2. Otherwise we have that $\nabla y\in H^{\tau-1}(\Omega)^2\hookrightarrow L^{r_p}(\Omega)^2$ for
\[\frac{1}{r_p} = 1-\frac{\tau}{2}\]
and, since $\tau > 1$, we can take $r_q = +\infty$. As before, we select
\[\frac{1}{s} = \min\left\{1-\frac{1}{\ps}-\frac{1}{r_p},1-\frac{1}{q}\right\}.\]
We have two possibilities now.

\medskip
\hspace{-\parindent}Step 3. If $\dfrac{1}{s} = 1-\dfrac{1}{q}$, then, applying \eqref{E3.4} and taking into account our choice of $q$, we have that $y\in H^{\hat\tau}(\Omega)$ with
\[\hat\tau = 1+\frac{2}{s} = 3-\frac{2}{q} = t,\]
and the proof is complete.

\medskip
\hspace{-\parindent}Step 4. Otherwise, $\dfrac{1}{s} = 1-\dfrac{1}{\ps}-\dfrac{1}{r_p}$ and we will have $y\in H^{\hat\tau}(\Omega)$ with
\[\hat\tau = 1+\frac{2}{s} = 1+2-\frac{2}{\ps}-(2-\tau) = \tau+1-\frac{2}{\ps},\]
and we have advanced a fixed amount $1-\dfrac{2}{\ps}$. If $\hat\tau\geq t$, the proof is complete.

\medskip
\hspace{-\parindent}Step 5. In other case, we can redefine $\tau = \hat\tau$ and go back to step 2.

\medskip
Every time we repeat the process, either we finish the proof or we increment the size of $\tau$ by the fixed amount $1-\frac{2}{\ps}$, so the proof will end in a finite number of steps.

\medskip
\hspace{-\parindent}(b)
From the Sobolev embedding theorem, we have that \[f\in L^r(\Omega)\hookrightarrow H^{2-t}(\Omega)'\text{ and }u\in \displaystyle\prod_{j=1}^m W^{1-1/r,r}(\Gamma_j)\hookrightarrow \displaystyle\prod_{j=1}^m H^{t-3/2}(\Gamma_j)\]
 for $t =\min\{2, 3-2/r\}$.
  The conditions imposed on $r$ imply that $1<t<1+\lambda$, $t\leq 2$, so we can apply Theorem \ref{T3.4}(a) to obtain $y\in H^t(\Omega)$ and we readily have that $y\in L^\infty(\Omega)$ and hence $a_0 y\in L^r(\Omega)$. Let us investigate the regularity of $b\cdot\nabla y$. We use again a boot-strap argument.

We have that $\nabla y\in H^{t-1}(\Omega)\hookrightarrow L^{\frac{2}{2-t}}(\Omega)$. Therefore $b\cdot\nabla y\in L^s(\Omega)$ where
\[\frac{1}{s} = \frac{1}{\ps}+\frac{2-t}{2}<1.\]
Applying Lemma \ref{L3.2}(b), we have that $y\in W^{2,\min\{s,r\}}(\Omega)$. If $s\geq r$, the proof is complete. Otherwise, we have that $\nabla y\in W^{1,s}(\Omega)\hookrightarrow L^{s^*}(\Omega)$, with
\[\frac{1}{s^*} = \frac{1}{s}-\frac{1}{2}.\]
Therefore
$b\cdot\nabla y\in L^{\hat s}(\Omega)$ where
\[\frac{1}{\hat s} = \frac{1}{\ps}+\frac{1}{s^*} = \frac{1}{\ps} + \frac{1}{s}-\frac{1}{2} = \frac{1}{s}-\left(\frac{1}{2}-\frac{1}{\ps}\right).\]
If $\frac{1}{\hat s}\leq \frac{1}{r}$, then the proof is complete. Otherwise, we can rename $s:=\hat s$ and repeat the argument subtracting at each step the positive constant $\dfrac{1}{2}-\dfrac{1}{\ps}$ until
$\dfrac{1}{\hat s}\leq \dfrac{1}{r}$.

\medskip
\hspace{-\parindent}(c) To obtain this result, we want to apply Lemma \ref{L3.2}(c), but the boundary datum in that result is in the space $\prod_{j = 1}^m V^{1-1/p,p}_{\vec\beta}(\Gamma_j)$, while the boundary datum in this result is in $\prod_{j = 1}^m W^{1-1/p,p}_{\vec\beta}(\Gamma_j)$. Taking into account Remark \ref{R3.3}, it is clear that for  $p<2$, the condition $\beta_j\geq 0$ implies that $\beta_j > 1-2/p$ and hence $W^{1-1/p,p}_{\vec\beta}(\Gamma_j) = V^{1-1/p,p}_{\vec\beta}(\Gamma_j)$ for all $j\in\{1,\ldots,m\}$.
If $p=2$, we define
\[u_s =\displaystyle\sum_{\beta_j>0} u\zeta_j,\]
where the $\zeta_j$ are the cutoff functions introduced in the proof of Lemma \ref{L3.2}(c). Taking into account again Remark \ref{R3.3} and noting that $u_s\equiv 0$ in a neighbourhood of the corners $\xS_j$ with $\beta_j=0$, it is readily deduced that $u_s\in \prod_{j=1}^m  V^{1-1/p,p}_{\vec\beta}(\Gamma_j)$. We also have that the function $u_r = u-u_s\in \prod_{j=1}^m W^{1-1/p,p}(\Gamma_j)$, because $u_r\equiv 0$ in a neighbourhood of the corners $\xS_j$ such that $\beta_j > 0$. In the same way we define
\[f_s = \displaystyle\sum_{\beta_j>0} f\zeta_j \in L^p_{\vec\beta}(\Omega)\text{ and } f_r= f-f_s\in L^p(\Omega),\]
and consider $y_s,y_r\in H^1(\Omega)$ such that
\[\langle \mathcal{A}y_r,z\rangle_\Omega = \langle f_r,z\rangle_\Omega + \langle  u_r,z\rangle_\Gamma,\text{ and }
\langle \mathcal{A}y_s,z\rangle_\Omega = \langle f_s,z\rangle_\Omega + \langle  u_s,z\rangle_\Gamma\quad \forall z\in H^1(\Omega),\]
so that $y =y_r+y_s$. As an application of Theorem \ref{T3.4}(b), $y_r\in W^{2,2}(\Omega)$, which is continuously embedded in $W^{2,2}_{\vec\beta}(\Omega)$ because $\beta_j\geq 0$ for all $j\in\{1,\ldots,m\}$.

Taking into account the above considerations, in the rest of the proof we assume that $\beta_j > 1-2/p$. If $p< 2$ then this holds, as discussed before. If $p=2$, we denote $u_s=u$ to treat both cases simultaneously, and hence we can use both that $u\in \prod_{j=1}^m W^{1-1/p,p}_{\vec\beta}(\Gamma_j)$, which is needed to have an embedding in a non-weighted Sobolev space, and $u\in \prod_{j=1}^m V^{1-1/p,p}_{\vec\beta}(\Gamma_j)$, which is needed to apply Lemma \ref{L3.2}(c).

From \cite[Lemma 2.29(ii)]{Pfefferer2014}, we deduce that $L^{p}_{\vec\beta}(\Omega)\hookrightarrow L^r(\Omega)$ for all $r<\frac{2}{{\beta_j}+2/p}\leq \frac{2}{2/p}= p$ for all $j\in\{1,\ldots,m\}$.
On the other hand, using the definition of the $\prod_{j = 1}^mW^{1-1/p,p}_{\vec\beta}(\Gamma_j)$-norm and \cite[Lemma 2.29(i)]{Pfefferer2014}, we have the embedding $\prod_{j = 1}^mW^{1-1/p,p}_{\vec\beta}(\Gamma_j)\hookrightarrow \prod_{j=1}^m W^{1-1/r,r}(\Gamma_j)$ for the same $r$ as above.
We notice at this point that
 the assumption $\beta_j<2-\frac{2}{p}$ implies that $\frac{2}{{\beta_j}+2/p}>1$, and
$2-\frac{2}{p}-\lambda_j <{\beta_j}$  implies $r < \dfrac{2}{2-\lambda_j}$ for all $j$, therefore we can choose some $r>1$ satisfying the assumptions of Theorem \ref{T3.4}(b) and we have that $a_0\in L^r(\Omega)$, $f\in L^r(\Omega)$, and $u\in\prod_{j=1}^m W^{1-1/r,r}(\Gamma_j)$. By Theorem \ref{T3.4}(b) we obtain $y\in W^{2,r}(\Omega)$ for
some $r>1$.

In particular, the result $y\in W^{2,r}(\Omega)$ implies $y\in L^\infty(\Omega)$, and hence $a_0 y \in L^{ p}_{{\vec\beta}}(\Omega)$.
We also have that $\nabla y\in W^{1,r}(\Omega)^2\hookrightarrow
L^{s_y}(\Omega)^2$ for $s_y = \dfrac{2r}{2-r}$ if $r<2$, any $s_y<+\infty$ if $r=2$ and $s_y=+\infty$ if $r>2$.
From this we deduce that $b\cdot \nabla y \in L^s(\Omega)$ for
\[\frac{1}{s} = \frac{1}{\ps}+\frac{1}{s_y}.\]
 Now we use that ${\vec\beta}\geq 0$ to deduce that $b\cdot\nabla y\in L^{s}_{{\vec\beta}}(\Omega)$ and hence $F = -b\cdot\nabla y -a_0 y \in L^{\min\{s,p\}}_{{\vec\beta}}(\Omega)$.
 By applying Lemma \ref{L3.2}(c), we have that $y\in W^{2,{\min\{s,p\}}}_{\vec\beta}(\Omega)$. If $s\geq  p$, the proof is complete.

Otherwise, in case $s < p \leq 2$, from Sobolev's embedding theorem, we have that $\nabla y\in W^{1,s}_{\vec\beta}(\Omega)\hookrightarrow L^{s_y}_{\vec\beta}(\Omega)$ for
 \[\frac{1}{s_y}=\frac{1}{s}-\frac{1}{2}=\frac{s-2}{s}\iff s_y = \frac{2s}{s-2}.\]
  Since $\vec\beta\geq \vec 0$, using that $b\in L^\ps(\Omega)$, we have that $b\cdot\nabla y \in L^{\hat s}_{\vec\beta}$, where
 \begin{equation}\label{E3.5}\frac{1}{\hat s} = \frac{1}{s_y}+\frac{1}{\ps} = \frac{1}{s}-\left(\frac{1}{2}-\frac{1}{\ps}\right). \end{equation}
 By applying Lemma \ref{L3.2}(c), we have that $y\in W^{2,\min\{p,\hat{s}\}}_{\vec\beta}(\Omega)$. If $\hat s\geq  2$, the proof is complete. Otherwise, we redefine $s:=\hat s$ and repeat the last step. Since at each iteration we subtract the positive constant $\dfrac{1}{2}-\dfrac{1}{\ps}$, the proof will end in a finite number of steps.
\end{proof}

We conjecture that the result of Theorem \ref{T3.4}(c) holds for $p\in(1,\hat p]$, but the proof is limited to $p\leq 2$.

Notice that the operator $\mathcal{A}^*$ is different from $\mathcal{A}$, and hence the results in Theorem \ref{T3.4} are not immediately applicable. For the adjoint state equation, we will need another assumption on $b\cdot n$, which is a result of the boundary term obtained due to integration by parts.
\begin{theorem}\label{T3.5}
Suppose assumptions \ref{A2.1} and \ref{A3.1} hold. Consider $f\in H^1(\Omega)'$ and $g\in H^{1/2}(\Gamma)'$ and  let $\varphi\in H^1(\Omega)$ be the unique solution of
\begin{equation}\label{E3.6}\langle \mathcal{A}^*\varphi,z\rangle_\Omega = \langle f,z\rangle_{\Omega}+\langle g,z\rangle_\Gamma\quad \forall z\in H^1(\Omega).\end{equation}

\medskip

\hspace{-\parindent}\textup{(a)} If $a_0,\, \nabla\cdot b\in L^q(\Omega)$, $b\cdot n\in L^{q_\Gamma}(\Gamma)\cap H^{t-3/2}(\Gamma)$, $f\in H^{2-t}(\Omega)'$, and $g\in \prod_{j=1}^m H^{t-3/2}(\Gamma_j)$  for $1<t<1+\lambda$, $t\leq 2$, $q = \dfrac{2}{3-t}$, and $q_\Gamma = \min\{2,1/(2-t)\}$, then $\varphi\in H^{t}(\Omega)$, and there exists a constant $C_{\mathcal{A}^*,t}>0$ such that
\[\Vert \varphi\Vert_{H^{t}(\Omega)} \leq C_{\mathcal{A}^*,t}\Big( \Vert f\Vert_{ H^{2-t}(\Omega)'} + \sum_{j=1}^m\Vert g\Vert_{H^{t-3/2}(\Gamma_j)}\Big).\]

\medskip

\hspace{-\parindent}\textup{(b)} If $a_0,\, \nabla\cdot b,\, f\in L^r(\Omega)$, and $g,b\cdot n\in \prod_{j=1}^m W^{1-1/r,r}(\Gamma_j)$ for some $r\in(1,\ps]$ satisfying $r<\dfrac{2}{2-\lambda}$ in case of $\lambda<2$, then $\varphi\in W^{2,r}(\Omega)$, and there exists a constant $C_{\mathcal{A}^*,r}>0$ such that
\[\Vert \varphi\Vert_{W^{2,r}(\Omega)}  \leq C_{\mathcal{A}^*,r}\Big(\Vert f\Vert_{ L^{r}(\Omega)} + \sum_{j=1}^m\Vert g\Vert_{W^{1-1/r,r}(\Gamma_j)}\Big).\]

\medskip

\hspace{-\parindent}\textup{(c)}   If $a_0,\, \nabla\cdot b,\, f\in L^{p}_{{\vec\beta}}(\Omega)$, and  $b\cdot n,\,g\in \prod_{j = 1}^m W^{1-1/p,p}_{\vec\beta}(\Gamma_j)$ for some $p\in(1,2]$ and some $\vec\beta$ such that $ 2-\frac{2}{p}-\lambda_j <\beta_j < 2-\frac{2}{p}$, $\beta_j \geq 0$, for all $j\in\{1,\ldots,m\}$,  then $\varphi\in W^{2,p}_{\vec\beta}(\Omega)$ and there exists a constant $C_{\mathcal{A}^*,\vec\beta,p}>0$ such that
\[\Vert \varphi\Vert_{W^{2,p}_{{\vec\beta}}(\Omega)}\leq C_{\mathcal{A}^*,\vec\beta,p}\Big(\Vert f\Vert_{ L^{2}_{\vec \beta}(\Omega)}+
\sum_{j=1}^m\Vert g\Vert  _{W^{1-1/p,p}_{\vec\beta}(\Gamma_j)}\Big).\]

\end{theorem}
\begin{proof}
The expression for $\langle \mathcal{A}^*\varphi,z\rangle_\Omega$ is derived in Lemma \ref{L2.5}. Using the product rule, we have that the function $\varphi$ satisfies
  \begin{align*}
    \int_\Omega \sum_{i,k = 1}^2 a_{ki}\partial_{x_i}\varphi \partial_{x_k} z \dx & -\int_\Omega (b\cdot\nabla \varphi) z \dx + \int_\Omega a_0 \varphi z\dx \\
    =  & \int_\Omega (\nabla \cdot b )\varphi z\dx -\int_\Gamma \varphi (b\cdot n) z\,\dx+ \langle f,z\rangle_\Omega
    +\langle g,z\rangle_\Gamma
  \end{align*}
  and we can apply Theorem \ref{T3.4} to this problem provided $\varphi\nabla\cdot b$ and $\varphi b\cdot n$ are in the appropriate spaces.

Notice that statement (a) for $t=2$ is the same than statement (b) for $r=2$. We will prove (a) for $t<2$, and refer to (b) for $t=2$.

{\em Step 1:} First, we prove $W^{1,\delta}(\Omega)$ regularity for some $\delta>2$.

Let us write the equation as
\[\left\{
\begin{array}{rcll}
A\varphi + \varphi &= &f+\varphi\nabla\cdot b + b\cdot\nabla\varphi- a_0\varphi + \varphi&\text{ in }\Omega\\
\partial_{n_A}\varphi &=& -b\cdot n \varphi+g&\text{ on }\Gamma.
\end{array}
\right.
\]
This is a Neumann problem posed on a Lipschitz domain. We will apply the regularity results in \cite{Geng2012}. To that end, we first investigate the existence of $r_f>2$ and $q_\Gamma > 1$ such that $f\in  W^{1,r_f'}(\Omega)'$ and $b\cdot n\in L^{q_\Gamma}(\Gamma)$.
In each of the three cases, we have:
\begin{itemize}
\item[(a)] $f\in H^{2-t}(\Omega)'\hookrightarrow W^{1,r_f'}(\Omega)'$ for $r_f=\dfrac{2}{2-t}>2$ since $1<t<2$. The exponent $q_\Gamma$ is given in the theorem.
\item[(b)] $f\in L^r(\Omega)\hookrightarrow W^{1,r_f'}(\Omega)'$ for $r_f = \dfrac{2r}{2-r} > 2$ if $1<r<2$ and all $r_f<+\infty$ if $r\geq 2$. In this case we take $q_\Gamma = r>1$.
\item[(c)] $f\in L^{p}_{\vec\beta}(\Omega)\subset L^r(\Omega)$ for $r<\dfrac{2}{\beta_j+\frac{2}{p}}$ for all $j\in\{1,\ldots,m\}$. The condition $\beta_j < 2-\dfrac{2}{p}$ implies that $\dfrac{2}{\beta_j+\frac{2}{p}} > 1$, so we can choose $r>1$ and $L^r(\Omega)\hookrightarrow W^{1,r_f'}(\Omega)'$ for $r_f = \dfrac{2r}{2-r} > 2$. Therefore $f\in W^{1,r_f'}(\Omega)'$ for all $2< r_f <\dfrac{2p}{2-(1-\beta_j) p}$.
     In this case we take $q_\Gamma = r>1$.
\end{itemize}
Note that also in each of the three cases we have different assumptions on $a_0$ and $\nabla\cdot b$, but in any case there exists $q_0>1$ such that $a_0, \nabla\cdot b\in L^{q_0}(\Omega)$.

Let us check that also $F = \varphi\nabla\cdot b + b\cdot\nabla\varphi - a_0\varphi +\varphi\in W^{1,r_\Omega'}(\Omega)'$ for some $r_\Omega>2$. To this end define $r_\varphi$, $s_\Omega$ and $r_\Omega$ by
\[\frac{1}{s_\Omega}=\frac{1}{r_\varphi}= \min\left\{
\frac{1}{2}\left(1-\frac{1}{q_0}\right), \frac{1}{2}-\frac{1}{p}\right\}\in (0,\tfrac{1}{2})\text{ and }
\frac{1}{r_\Omega} = \frac{1}{2}-\frac{1}{s_\Omega}\in (0,\tfrac{1}{2})\]
such that
\[\frac{1}{r_\varphi}+\frac{1}{q_0}+\frac{1}{s_\Omega}\leq 1\text{ and }\frac{1}{p}+\frac{1}{2}+\frac{1}{s_\Omega}\leq 1\]
and $W^{1,r_\Omega'}(\Omega)\hookrightarrow L^{s_\Omega}(\Omega)$.
Using Lemma \ref{L2.6}, we have that $\varphi\in H^1(\Omega)\hookrightarrow L^{r_\varphi}(\Omega)$. Denote $S = \{z\in W^{1,r_\Omega'}(\Omega): \Vert z\Vert_{W^{1,r_\Omega'}(\Omega)} = 1\}$. Then
\begin{align*}
  \Vert  &F\Vert_{W^{1,r_\Omega'}(\Omega)'} =   \sup_{z\in S}\int_\Omega \left(\varphi\nabla\cdot b + b\cdot\nabla\varphi + a_0\varphi -\varphi\right) z\dx\\
  &\leq  C\sup_{z\in S} \big(\Vert \varphi\Vert_{L^{r_\varphi}(\Omega)} \Vert 1+a_0+\nabla\cdot b\Vert_{L^{q_0}(\Omega)} + \Vert b\Vert_{L^p(\Omega)}\Vert \nabla\varphi\Vert_{L^2(\Omega)} \big)\Vert z\Vert_{L^{s_\Omega}(\Omega)}\\
  &\leq  C_{r_\Omega}\sup_{z\in S}\big(\Vert \varphi\Vert_{L^{r_\varphi}(\Omega)} (\Vert 1+a_0+\nabla\cdot b\Vert_{L^{q_0}(\Omega)} )+ \Vert b\Vert_{L^p(\Omega)}\Vert \nabla\varphi\Vert_{L^2(\Omega)} \big)\Vert z\Vert_{W^{1,r_\Omega'}(\Omega)}\\
   & =  C_{r_\Omega}\big(\Vert \varphi\Vert_{L^{r_\varphi}(\Omega)} (\Vert 1+a_0+\nabla\cdot b\Vert_{L^{q_0}(\Omega)} )+ \Vert b\Vert_{L^p(\Omega)}\Vert \nabla\varphi\Vert_{L^2(\Omega)} \big).
\end{align*}

On the boundary, we want to check that $b\cdot n \varphi\in W^{-1/r_\Gamma,r_\Gamma}(\Gamma) = W^{1/r_\Gamma,r_\Gamma'}(\Gamma)'$ for some $r_\Gamma>2$. To this end, define $\hat{r}_\varphi$, $s_\Gamma$ and $r_\Gamma$ by
\[\frac{1}{s_\Gamma}=\frac{1}{\hat{r}_\varphi} = \frac{1}{2}\left(1-\frac{1}{q_\Gamma}\right) \in(0,\tfrac{1}{2})\text{ and }\frac{1}{r_\Gamma} = \frac{1}{2}\left(1-\frac{1}{s_\Gamma}\right)\in (0,\tfrac{1}{2})\]
such that
\[\frac{1}{\hat{r}_\varphi}+\frac{1}{q_\Gamma}+\frac{1}{s_\Gamma} = 1\]
and $W^{1/r_\Gamma,r_\Gamma'}(\Gamma)\hookrightarrow L^{s_\Gamma}(\Gamma)$. From Lemma 2.3 and the trace theorem, we have that
$\varphi\in H^{1/2}(\Gamma)\hookrightarrow L^{\hat{r}_\varphi}(\Gamma)$. Denote $S = \{z\in W^{1/r_\Gamma,r_\Gamma'}(\Gamma):\ \Vert z\Vert_{W^{1/r_\Gamma,r_\Gamma'}(\Gamma)} =1\}$. Then
\begin{align*}
  \Vert b&\cdot n \varphi\Vert_{W^{-1/r_\Gamma,r_\Gamma}(\Gamma)} = \sup_{z\in S}\int_\Gamma b\cdot n \varphi z\,\dx
  \leq \sup_{z\in S} \Vert b\cdot n\Vert_{L^{q_\Gamma}(\Gamma)} \Vert \varphi\Vert_{L^{\hat{r}_\varphi}(\Gamma)} \Vert z\Vert_{L^{s_\Gamma}(\Gamma)}\\
  &\leq  C_{r_\Gamma} \sup_{z\in S} \Vert b\cdot n\Vert_{L^{q_\Gamma}(\Gamma)} \Vert \varphi\Vert_{L^{\hat{r}_\varphi}(\Gamma)} \Vert z\Vert_{W^{1/r_\Gamma,r_\Gamma'}(\Gamma)}(\Gamma)  = C_{r_\Gamma} \Vert \varphi\Vert_{L^{\hat{r}_\varphi}(\Gamma)}\Vert b\cdot n\Vert_{L^{q_\Gamma}(\Gamma)}.
\end{align*}
Noting that for a general Lipschitz domain the $W^{1,\delta}(\Omega)$ regularity is limited to $\delta\leq 4$, see \cite{Geng2012}, from the previous estimates, we can deduce that, for $\delta = \min\{4,r_f,r_\Omega,r_\Gamma\}>2$, $\varphi\in W^{1,\delta}(\Omega)$.

{\em Step 2:} Let us check that $\varphi\nabla\cdot b$ and $\varphi b\cdot n$ satisfy the regularity assumptions for the source and the Neumann data respectively of the different cases of Theorem \ref{T3.4}.

\hspace{-\parindent}(a)  On one hand $\varphi\in W^{1,\delta}(\Omega)\hookrightarrow L^\infty(\Omega)$ and the assumption $\nabla\cdot b\in L^q(\Omega)$ imply $\varphi\nabla\cdot b\in L^q(\Omega)\hookrightarrow H^{2-t}(\Omega)'$, by the definition of $q$.
  On the boundary, by the trace theorem $\varphi\in W^{1-1/\delta,\delta}(\Gamma)$. If $1<t\leq 3/2$, then we use that $W^{1-1/\delta,\delta}(\Gamma)\hookrightarrow L^\infty(\Gamma)$ to conclude that $\varphi b\cdot n\in L^{q_\Gamma}(\Gamma)\hookrightarrow H^{t-3/2}(\Gamma)$. The last inclusion follows by duality and the Sobolev imbedding $H^{3/2-t}(\Gamma)\hookrightarrow L^{\frac{1}{t-1}}(\Gamma)$. If $3/2<t<2$, we use that
  $W^{1-1/\delta,\delta}(\Gamma)\hookrightarrow H^{s_1}(\Gamma)$ for $s_1=1-1/\delta> 1/2$. Since we are assuming that $b\cdot n\in H^{s_2}(\Gamma)$ with $s_2 = t-3/2\in (0,1/2)$, from the trace theorem and the multiplication theorem \cite[Theorem 7.4]{BehzadanHolst2021}, we have that $\varphi b\cdot n\in H^{t-3/2}(\Gamma)$.

The result follows from Theorem \ref{T3.4}(a).

\medskip

\hspace{-\parindent}(b) Using again that $\varphi\in L^\infty(\Omega)$, we readily deduce that $\varphi\nabla\cdot b\in L^r(\Omega)$.
Let us check that $\varphi b\cdot n\in\prod_{j = 1}^m W^{1-1/r,r}(\Gamma_j)$.

For all $j\in\{1,\ldots,m\}$, by the trace theorem and the assumption on $b\cdot n$ we deduce the existence of $B_j\in W^{1,r}(\Omega)$ such that the trace of $B_j$ on $\Gamma_j$ is $b\cdot n$.

Suppose first that $r\leq 2$. Then, a straightforward application of the multiplication Lemma \ref{L3.6} below (in the case $\beta_j=0$) yields $\varphi B_j\in W^{1,r}(\Omega)$, and hence, its trace on $\Gamma_j$ belongs to $W^{1-1/r,r}(\Gamma_j)$. Therefore, $\varphi b\cdot n\in\prod_{j = 1}^m W^{1-1/r,r}(\Gamma_j)$ and the result follows from Theorem \ref{T3.4}(b).

If $r > 2$, from the previous paragraph we have that $\varphi\in W^{2,2}(\Omega)\hookrightarrow W^{1,\delta}(\Omega)$ for all $\delta < +\infty$. Repeating the previous argument, we obtain the desired result.

\medskip

\hspace{-\parindent}(c) Since $\varphi\in L^\infty(\Omega)$ and $\nabla\cdot b\in L^p_{\vec\beta}(\Omega)$, we have that $\varphi\nabla\cdot b\in L^p_{\vec\beta}(\Omega)$. Next, we show that $\varphi b\cdot n\in\prod_{j = 1}^m W^{1-1/p,p}_{\vec\beta}(\Gamma_j)$.

For all $j\in\{1,\ldots,m\}$, by the trace theorem and the assumption on $b\cdot n$ we deduce the existence of $B_j\in W^{1,p}_{\vec\beta}(\Omega)$ such that the trace of $B_j$ on $\Gamma_j$ is $b\cdot n$.

Since $p\leq 2<\delta$, a straightforward application of the multiplication Lemma \ref{L3.6} below yields $\varphi B_j\in W^{1,p}_{\vec\beta}(\Omega)$, and hence, its trace on $\Gamma_j$ belongs to $W^{1-1/p,p}_{\vec\beta}(\Gamma_j)$. Therefore, $\varphi b\cdot n\in\prod_{j = 1}^m W^{1-1/p,p}_{\vec\beta}(\Gamma_j)$ and the result follows from Theorem \ref{T3.4}(c).
\end{proof}

It remains to prove the multiplication theorem used in the proofs of cases (b) and (c) in Theorem  \ref{T3.5}.

\begin{lemma}[A multiplication theorem in weighted Sobolev spaces]\label{L3.6}
Let $1<q < +\infty$.
Consider $\varphi\in W^{1,\delta}(\Omega)$ for some $\delta >\max\{2, q \}$ and $\psi\in W^{1,q }_{{\vec\beta}}(\Omega)$ for some $\vec\beta$ such that $2-\frac{2}{q }-\lambda_j < \beta_j <2-\frac{2}{q }$, $\beta_j\geq 0$
for all $j\in\{1,\ldots,m\}$. Then $\psi\varphi \in W^{1,q }_{{\vec\beta}}(\Omega)$.
\end{lemma}
\begin{proof}
  Since $\delta > 2$, $\varphi\in L^\infty(\Omega)$. Also it is clear that $\psi\in L^{q}_{\vec\beta}(\Omega)$, and hence $\psi\varphi\in L^{q}_{\vec\beta}(\Omega)$.

  Let us check that also $\vert\nabla (\psi\varphi)\vert\in L^{q}_{\vec\beta}(\Omega)$. We write $\nabla (\psi\varphi) = \varphi\nabla \psi + \psi\nabla\varphi$. Using again that $\varphi\in L^\infty(\Omega)$ it is immediate to deduce that  $\vert\nabla \psi\vert\in L^{q}_{\vec\beta}(\Omega)$ implies that  $\vert\varphi\nabla \psi\vert\in L^{q}_{\vec\beta}(\Omega)$.

  Checking that  the term $\psi\vert\nabla\varphi\vert\in L^{q}_{\vec\beta}(\Omega)$ is more involved. By localizing the problem at corner $x_j$, and applying H\"older's inequality we obtain
  \[\int_{{\Omega_{R_j}}} (r^{{\beta_j}} \psi \vert\nabla \varphi\vert)^q\dx \leq \Vert r^{{\beta_j}} \psi \Vert_{L^{\frac{{q}\delta}{\delta-{q}}}({\Omega_{R_j}})}^{q}\Vert \nabla \varphi\Vert_{L^\delta({\Omega_{R_j}})}^{q},\]
  and therefore it is sufficient to prove that $r^{\beta_j}\psi\in L^{\frac{q \delta}{\delta-{q}}}(\Omega_{R_j})$.
  Let us introduce $1\leq q_\delta< q $ and $2\leq q_\delta^*<+\infty$ such that
  \[\frac{1}{q_\delta^*}=\min\left\{\frac{1}{2},\frac{1}{q}-\frac{1}{\delta}\right\}\text{ and }  \frac{1}{q_\delta}=\min\left\{1,\frac{1}{q}+\frac{1}{2}-\frac{1}{\delta}\right\}=\frac{1}{q_\delta^*}+\frac{1}{2}\]
  so that  $q_\delta^* \geq  \frac{q \delta}{\delta-q}$, and $W^{1,q_\delta}(\Omega_{R_j})\hookrightarrow L^{q_\delta^*}(\Omega_{R_j})\hookrightarrow L^{\frac{q \delta}{\delta-{q}}}(\Omega_{R_j})$.
  We are going to prove that $r^{\beta_j}\psi\in W^{1,q_{\delta}}(\Omega_{R_j})$.

  First of all we notice that $\nabla(r^{{\beta_j}} \psi) = r^{{\beta_j}} \nabla \psi + \psi\nabla r^{{\beta_j}}$.
  By definition of $W^{1,q}_{\vec\beta}(\Omega)$, we have that $ r^{\beta_j} \vert\nabla \psi\vert \in L^{q}(\Omega_{R_j})\hookrightarrow L^{q_\delta}(\Omega_{R_j})$.

  For the second term we notice that $\vert\psi\nabla r^{\beta_j}\vert \sim r^{{{\beta_j}}-1} \psi$.
  Since $1-2/q_\delta=\max\{-1,\frac{2}{\delta}-\frac{2}{q}\} <0\leq\beta_j$, we have that
  $W^{1,q_{\delta}}_{\vec\beta}(\Omega_{R_j})\hookrightarrow L^{q_{\delta}}_{\vec\beta-1}(\Omega_{R_j})$; see e.g. \cite[Lemma 2.29(i)]{Pfefferer2014}.
    We deduce that $\psi\in W^{1,q}_{\vec\beta}(\Omega_{R_j}) \hookrightarrow W^{1,q_\delta}_{\vec\beta}(\Omega_{R_j})\hookrightarrow L^{q_{\delta}}_{\vec\beta-1}(\Omega_{R_j})$. This means that $r^{\beta_j-1} \psi\in L^{q_{\delta}}(\Omega_{R_j})$, and we gather that $\vert\psi\nabla r^{{\beta_j}}\vert\in L^{q_{\delta}}(\Omega_{R_j})$.

    Therefore $\nabla(r^{\beta_j} \psi)\in L^{q_{\delta}}(\Omega_{R_j})$, so we have that $r^{\beta_j} \psi\in W^{1,q_{\delta}}({\Omega_{R_j}})$.

  Using this,
  we conclude that $\psi\vert\nabla\varphi\vert\in L^{q}_{{\vec\beta}}(\Omega)$ and consequently $\vert \nabla (\psi\varphi)\vert \in L^{q}_{{\vec\beta}}(\Omega)$, which leads to the desired result.
\end{proof}

\section{Discretization}\label{S4}
Consider a family of  regular triangulations $\{\mathcal{T}_h\}$  graded with mesh grading parameters $\mu_j\in(0,1]$, $j\in\{1,\ldots,m\}$ in the sense of \cite[Section 3.1]{ASW1996}, see also \cite{AMPR2019}.
As usual, $Y_h\subset H^1(\Omega)\cap C(\bar\Omega)$ is the space of continuous piecewise linear functions.

\begin{lemma}\label{L4.1} There exists a constant $c_{\vec\mu}>0$ such that
  \[\Vert \psi- I_h\psi\Vert_{H^1(\Omega)}\leq c_{\vec\mu} h^s \Vert \psi\Vert_{W^{2,2}_{\vec\beta}(\Omega)}\ \forall  \psi\in W^{2,2}_{\vec\beta}(\Omega),\]
    where $I_h $ is the Lagrange interpolation operator, the vector  $\vec\beta$ satisfies that $1-\lambda_j< \beta_j < 1$ and $\beta_j\geq 0$ for all $j\in\{1,\ldots,m\}$
    and the exponent $s$ satisfies that $s\leq 1$ and $s <\frac{\lambda_j}{\mu_j}$ for all $j\in \{1,\ldots,m\}$.
\end{lemma}
\begin{proof}
  The case $\mu_j=1$ (quasi-uniform mesh) is classical. For $\mu_j<\lambda_j$ see \cite[Lemma 4.1]{APR2012}. The case $\lambda_j\leq\mu_j<1$ can be proved with the same techniques and the additional idea that $h_T\sim  h^s r_T^{1-s\mu_j}$, $1-s\mu_j=\beta_j > 1-\lambda_j$; see equation (3.14) in \cite[Theorem 3.2]{ASW1996}, where it was used for a Dirichlet problem.
\end{proof}

Define the bilinear form $a(y,z)= \langle\mathcal{A}y,z\rangle_\Omega$. For a datum $u\in H^{1/2}(\Gamma)'$, the discrete state equation reads
\begin{equation}\label{E4.1}a(y_h,z_h)=\langle u, z_h\rangle_\Gamma\ \forall z_h\in Y_h.\end{equation}
Existence and uniqueness of the solution of this equation is not immediate since $a(\cdot,\cdot)$ is not coercive over $Y_h$.

\begin{theorem}\label{T4.2}
  There exists $h_0>0$ that depends on $A$, $b$, $a_0$, $\Omega$ and the mesh grading parameter $\vec\mu$, such that the system \eqref{E4.1}
  has a unique solution for every $h\leq h_0$ and every $u\in H^{1/2}(\Gamma)'$. Further, there exists a constant $K_0$ that depends on $A$, $b$, $a_0$, $\Omega$ and is independent of $\vec\mu$ and $h$ such that
  \begin{equation}\label{E4.2}
  \Vert y_h\Vert_{H^1(\Omega)}\leq K_0\Vert \mathcal{A}^{-1} u\Vert_{H^1(\Omega)} \ \forall h\leq h_0.
  \end{equation}
\end{theorem}
The scheme of the proof is similar to that of \cite[Lemma 3.1]{CMR2021} for distributed control problems with homogeneous Dirichlet boundary conditions, but that proof is done for quasi-uniform meshes and uses this fact explicitly; see equations (3.8) and (3.9) in \cite{CMR2021}. Since the mesh grading and the boundary terms imply some extra technicalities, we include a complete proof for the convenience of the reader.
\begin{proof}
  Due to the linearity of the system, to show existence it is sufficient to prove uniqueness of solution in the case $u=0$. Suppose $y_h\in Y_h$ satisfies
  \begin{equation}\label{E4.3}a(y_h,z_h) = 0\ \forall z_h\in Y_h.\end{equation}
  Taking $z_h=y_h$ and using G\r{a}rding's inequality established in Lemma \ref{L2.3}, we have that
  \[0=a(y_h,y_h)= \langle \mathcal{A} y_h,y_h\rangle_\Omega \geq
   \frac{\Lambda}{8 C_E^2}\Vert y_h\Vert^2_{H^1(\Omega)} - C_{\Lambda,E,b}\Vert y_h\Vert^2_{L^2(\Omega)}.
   \]
   Therefore
   \begin{equation}\label{E4.4}
   \Vert y_h\Vert_{H^1(\Omega)}\leq 2C_E\sqrt{\frac{ 2 C_{\Lambda,E,b}}{\Lambda}}\Vert y_h\Vert_{L^2(\Omega)}.\end{equation}
   Since $y_h\in L^2(\Omega)\subset L^2_{\vec\beta}(\Omega)$ for all $\vec\beta\geq \vec 0$ such that $1-\lambda_j < \beta_j $ for all $j\in\{1,\ldots,m\}$, from Theorem \ref{T3.5}(c), we have that there exists a unique $\psi\in W^{2,2}_{\vec\beta}(\Omega)$ such that
   \begin{equation}\label{E4.5}a(z,\psi) =
   \int_\Omega y_h z\mathrm{d}x\ \forall z\in H^1(\Omega)\end{equation}
   and there exists a constant $C_{\mathcal{A}^*,\vec\beta}$ such that
   \[\Vert \psi\Vert_{W^{2,2}_{\vec\beta}(\Omega)}\leq C_{\mathcal{A}^*,\vec\beta}\Vert y_h\Vert_{L^2_{\vec\beta}(\Omega)}.\]
   Let us denote $\hat \psi_h\in Y_h$ the Ritz-Galerkin projection of $\psi$ onto $Y_h$ in the sense of $H^1(\Omega)$, i.e., $\hat \psi_h$ is the unique solution of
   \[\int_\Omega (\nabla \hat \psi_h\nabla z_h + \hat\psi_h z_h )\mathrm{d}x =
   \int_\Omega (\nabla \psi\nabla z + \psi z )\mathrm{d}x .\]
   From \cite[Eq. (4.2)]{APR2015}, Theorem \ref{T3.5}(c), and the embedding $L^2(\Omega)\hookrightarrow L^2_{\vec\beta}(\Omega)$, with embedding constant~$1$ due to the choice $R_j\leq 1$, we have that there exists a constant $\hat c_{\vec\mu}$ such that
   \begin{equation}
   \label{E4.6}
   \Vert \psi-\hat\psi_h\Vert_{H^1(\Omega)}\leq \hat c_{\vec\mu} h^{s}\Vert \psi\Vert_{W^{2,2}_\beta(\Omega)}
   \leq \hat c_{\vec\mu} C_{\mathcal{A}^*,\beta} h^s \Vert y_h\Vert_{L^2_{\vec\beta}(\Omega)} \leq \hat c_{\vec\mu} C_{\mathcal{A}^*,\beta} h^s \Vert y_h\Vert_{L^2(\Omega)},
   \end{equation}
where $s\leq 1$ and $s <\frac{\lambda_j}{\mu_j}$ for all $j\in \{1,\ldots,m\}$; see Lemma \ref{L4.1}.
   Taking $z=y_h$ in the adjoint equation \eqref{E4.5}, and $z_h=\hat\psi_h$ in the homogeneous discrete equation \eqref{E4.3}, we deduce
   \begin{align*}
   \Vert y_h\Vert^2_{L^2(\Omega)} &= a(y_h,\psi) = a(y_h,\psi-\hat\psi_h) \leq \Vert \mathcal{A}\Vert \Vert y_h\Vert_{H^1(\Omega)}\Vert \psi-\hat\psi_h\Vert_{H^1(\Omega)}\\
   &\leq \hat c_{\vec\mu} C_{\mathcal{A}^*,\vec\beta} \Vert \mathcal{A}\Vert \Vert y_h\Vert_{H^1(\Omega)} h^s \Vert y_h\Vert_{L^2(\Omega)}.
   \end{align*}
   Along the proof we will denote $\Vert \mathcal{A}\Vert  = \Vert \mathcal{A}\Vert_{\mathcal{L}(H^1(\Omega),H^1(\Omega)')}$.
   Choosing $h_0$ such that
   \begin{equation}\label{E4.7} \hat c_{\vec\mu} C_{\mathcal{A}^*,\vec\beta} \Vert \mathcal{A} \Vert  h_0^s=\frac{1}{2}\frac{1}{2C_E} \sqrt{ \frac{\Lambda}{ 2 C_{\Lambda,E,b}}},\end{equation}
   we have that, for all $h \leq  h_0$
   \[ \Vert y_h\Vert_{L^2(\Omega)}\leq \frac{1}{2}\frac{1}{2C_E} \sqrt{ \frac{\Lambda}{ 2 C_{\Lambda,E,b}}} \Vert y_h\Vert_{H^1(\Omega)}.\]
   Using this and estimate \eqref{E4.4}, we deduce that \[\Vert y_h\Vert_{H^1(\Omega)}\leq \frac{1}{2}\Vert y_h\Vert_{H^1(\Omega)}\ \forall h \leq  h_0,\]
   and hence $y_h=0$.

   Take now $u\in H^{1/2}(\Gamma)'$ and denote $y=\mathcal{A}^{-1}u$. For $h\leq h_0$, let $y_h$ be the solution of \eqref{E4.1}. Taking $z=y_h$ in the adjoint equation \eqref{E4.5}, and $z_h=\hat\psi_h$ in the discrete equation \eqref{E4.1}, we deduce
   \begin{align*}
   \Vert y_h\Vert^2_{L^2(\Omega)} &= a(y_h,\psi) = a(y_h,\psi-\hat\psi_h)+
   \langle u,\hat\psi_h\rangle_\Gamma
    =  a(y_h,\psi-\hat\psi_h)+a(y,\hat\psi_h)\\
    &\leq  \Vert \mathcal{A}\Vert \left( \Vert y_h\Vert_{H^1(\Omega)} \Vert \psi-\hat\psi_h\Vert_{H^1(\Omega)} + \Vert y\Vert_{H^1(\Omega)}\Vert \hat \psi_h\Vert_{H^1(\Omega)}\right)\\
    &\leq  \hat c_{\vec\mu} C_{\mathcal{A}^*,\vec\beta} \Vert \mathcal{A}\Vert \Vert y_h\Vert_{H^1(\Omega)} h^s \Vert y_h\Vert_{L^2(\Omega)} + c_{\vec \beta} C_{\mathcal{A}^*,\vec\beta} \Vert \mathcal{A}\Vert  \Vert y\Vert_{H^1(\Omega)}\Vert y_h\Vert_{L^2(\Omega)},
    \end{align*}
    where we have used that $\Vert \hat\psi_h\Vert_{H^1(\Omega)}\leq \Vert \psi\Vert_{H^1(\Omega)}\leq \hat c_{\vec \beta}\Vert \psi\Vert_{W^{2,2}_{\vec\beta}(\Omega)}\leq c_{\vec \beta} C_{\mathcal{A}^*,\vec\beta}\Vert y_h\Vert_{L^2(\Omega)}$; see \cite[Lemma 2.29(i)]{Pfefferer2014} for the embedding $W^{2,2}_{\vec\beta}(\Omega)\hookrightarrow H^1(\Omega)$.
    Now, using that $h\leq h_0$ and \eqref{E4.7} we have
    \begin{align*}
      \Vert y_h\Vert_{L^2(\Omega)} &\leq   \frac{1}{2}\frac{1}{2C_E} \sqrt{ \frac{\Lambda}{ 2 C_{\Lambda,E,b}}} \Vert y_h\Vert_{H^1(\Omega)}
        +c_{\vec \beta} C_{\mathcal{A}^*,\vec\beta} \Vert \mathcal{A}\Vert  \Vert y\Vert_{H^1(\Omega)},
    \end{align*}
    and applying Young's inequality we deduce
    \begin{equation}\label{E4.8}
      \Vert y_h\Vert_{L^2(\Omega)}^2 \leq   \frac{1}{16}\frac{1}{C_E^2}  \frac{\Lambda}{  C_{\Lambda,E,b}} \Vert y_h\Vert_{H^1(\Omega)}^2
        +2 c_{\vec \beta}^2C_{\mathcal{A}^* ,\vec\beta}^2 \Vert \mathcal{A}\Vert^2 \Vert y\Vert_{H^1(\Omega)}^2.
    \end{equation}
    Using G\r{a}rding's inequality, the discrete equation \eqref{E4.1} and $y=\mathcal{A}^{-1}u$, we infer
    \begin{align}
    \frac{\Lambda}{8 C_E^2}\Vert y_h\Vert^2_{H^1(\Omega)} &- C_{\Lambda,E,b}\Vert y_h\Vert^2_{L^2(\Omega)}  \leq a(y_h,y_h)\nonumber \\
    & =
    \langle u,y_h\rangle_\Gamma
    = a(y,y_h) \leq \Vert \mathcal{A}\Vert \Vert y\Vert_{H^1(\Omega)} \Vert y_h\Vert_{H^1(\Omega)}.\label{E4.9}
    \end{align}
    Multiplying \eqref{E4.8} by $C_{\Lambda,E,b}$ and using the resulting inequality in \eqref{E4.9}, we obtain
    \begin{align*}
    \frac{\Lambda}{16C_E^2} &\Vert y_h\Vert_{H^1(\Omega)}^2 \leq 2 c_{\vec \beta}^2C_{\mathcal{A}^*,\vec \beta}^2 \Vert \mathcal{A}\Vert^2 \Vert y\Vert_{H^1(\Omega)}^2  +
    \Vert \mathcal{A}\Vert \Vert y\Vert_{H^1(\Omega)} \Vert y_h\Vert_{H^1(\Omega)}\\
    &\leq 2 c_{\vec \beta}^2 C_{\mathcal{A}^*,\vec \beta}^2 \Vert \mathcal{A}\Vert^2 \Vert y\Vert_{H^1(\Omega)}^2  +\frac{8 C_E^2}{\Lambda}\Vert \mathcal{A}\Vert^2 \Vert y\Vert_{H^1(\Omega)}^2 + \frac{\Lambda}{32C_E^2} \Vert y_h\Vert_{H^1(\Omega)}^2,
    \end{align*}
    where in the second step we have used Young's inequality.
    Gathering the terms with $\Vert y_h\Vert_{H^1(\Omega)}^2$ and taking the square root, we finally obtain:
    \begin{align*}\Vert y_h\Vert_{H^1(\Omega)} &\leq \frac{C_E}{4\sqrt{2\Lambda}} \Vert \mathcal{A}\Vert \left(2c_{\vec \beta}^2 C_{\mathcal{A}^*,\vec \beta}^2 + \frac{8 C_E^2}{\Lambda}\right)^{1/2} \Vert \mathcal{A}^{-1}u\Vert_{H^1(\Omega)}.
    \end{align*}
    Notice that the constant depends on $\vec\beta$, which is itself limited by the value of $\vec\lambda$, and hence the constant will finally depend on $\vec\lambda$.
\end{proof}

\begin{theorem}\label{T4.3}
There exists $h_0^*>0$ that depends on $A$, $b$, $a_0$, $\Omega$ and the mesh grading parameter $\vec\mu$, such that the discrete adjoint problem
  \begin{equation}\label{E4.10}a(z_h,\varphi_h)=\langle y,z_h\rangle_\Omega\ \forall z_h\in Y_h\end{equation}
  has a unique solution for every $y\in H^1(\Omega)'$ and every $0< h \leq h_0^*$. Further, there  exists a constant $K_0^*$ that depends on $A$, $b$, $a_0$, $\Omega$ and is independent of $\vec\mu$ and $h$ such that
  \begin{equation}\label{E4.11}\Vert \varphi_h\Vert_{H^1(\Omega)}\leq K_0^*\Vert (\mathcal{A}^{*})^{-1}y\Vert_{H^1(\Omega)}\ \forall h<h_0^*.
  \end{equation}
\end{theorem}
\begin{proof}
    Existence and uniqueness of solution of the discrete adjoint equation \eqref{E4.10} follows for all $0<h<h_0$ due to the finite-dimensional character of the problem.
    To get the estimate \eqref{E4.11}, we follow the steps of the proof of Theorem \ref{T4.2}. Notice that in this case, the value of $h_0^*$, which is used explicitly in the proof, may be different from the value of $h_0$ provided in \eqref{E4.7}.
\end{proof}
The following estimate is an immediate consequence of the previous results, Lemma \ref{L2.6}, Corollary \ref{C2.7} and the trace theorem.
\begin{corollary}\label{C4.4}
Let $\bar h=\min\{h_0,h_0^*\}$ with $h_0$ from Theorem \ref{T4.2} and $h_0^*$ from Theorem \ref{T4.3}. For $u\in L^2(\Gamma)$ let $y_h\in Y_h$ be the unique solution of \eqref{E4.1}.
There exists a constant $c_2>0$ that depends on the data of the problem, but not on the mesh grading parameters $\vec\mu$ or on $h$, such that, for all $h<\bar h$
\begin{align}
  \Vert y_h\Vert_{L^2(\Gamma)} & \leq c_2\Vert u\Vert_{L^2(\Gamma)}. \label{E4.12}
\end{align}
\end{corollary}
\begin{proof}

   Let us denote $C_{\mathrm{TR}}$ the norm of the trace operator from $H^1(\Omega)$ to $L^2(\Gamma)$. We use Theorem \ref{T4.2}, Lemma \ref{L2.6}, and the fact that $u$ can be seen as an element of $H^1(\Omega)'$ and $\Vert u\Vert_{H^1(\Omega)'} \leq C_{\mathrm{TR}} \Vert u\Vert_{L^2(\Gamma)}$, cf. \eqref{E2.4} and \eqref{E2.5}. A straightforward estimation shows that
  \begin{align*}
    \Vert y_h\Vert_{L^2(\Gamma)} & \leq C_{\mathrm{TR}}\Vert y_h\Vert_{H^1(\Omega)} \leq C_{\mathrm{TR}} K_0\Vert \mathcal{A}^{-1}u\Vert_{H^1(\Omega)} \leq C_{\mathrm{TR}} K_0 \Vert \mathcal{A}^{-1}\Vert  C_{\mathrm{TR}} \Vert u\Vert_{L^2(\Gamma)},
  \end{align*}
  where $\Vert \mathcal{A}^{-1}\Vert $ denotes the norm in $\mathcal{L}(H^1(\Omega)',H^1(\Omega))$.
  The result follows for $c_2 = C_{\mathrm{TR}}^2 K_0 \Vert \mathcal{A}^{-1}\Vert $.
\end{proof}
\begin{theorem}[Error estimates in the domain] \label{T4.5}
For $0< h<\bar h$, where $\bar h$ is defined in Corollary \ref{C4.4}, and $u\in H^{1/2}(\Gamma)'$, let $y_h\in Y_h$ be the solution of \eqref{E4.1} and $y\in H^1(\Omega)$ be the solution of \eqref{E3.1} for $f=0$. There exists $C>0$ that depends on $A$, $b$, $a_0$, $\Omega$ but is independent of $h$ such that
  \begin{equation}\label{E4.13}\Vert y-y_h\Vert_{L^2(\Omega)} \leq C h^s \Vert u\Vert_{H^{1/2}(\Gamma)'}.\end{equation}
If further
$u\in \prod_{j=1}^m W^{1/2,2}_{\vec\beta}(\Gamma_j)$, where $1-\lambda_j < \beta_j < 1$ and $\beta_j\geq 0$ for all $j\in\{1,\ldots,m\}$, there exists $C>0$ that depends on $A$, $b$, $a_0$, $\Omega$, $\vec\beta$, and the mesh grading parameter $\vec\mu$, but is independent of $h$ and $u$ such that
  \begin{equation}\label{E4.14}\Vert y-y_h\Vert_{L^2(\Omega)}+ h^s \Vert y-y_h\Vert_{H^1(\Omega)} \leq C h^{2s} \Vert y\Vert_{W^{2,2}_{\vec\beta}(\Omega)} \leq C h^{2s} \sum_{j=1}^m\Vert u\Vert_{ W^{1/2,2}_{\vec\beta}(\Gamma_j)}\end{equation}
 for all $s\leq 1$ and $s <\frac{\lambda_j}{\mu_j}$ for all $j\in \{1,\ldots,m\}$.

Furthermore, for all $f\in L^{2}_{\vec\beta}(\Omega)$ and $g\in \prod_{j=1}^m W^{1/2,2}_{\vec\beta}(\Gamma)$, let $\varphi\in W^{2,2}_{\vec\beta}(\Omega)$ be the solution of \eqref{E3.6} and $\varphi_h$ be the unique solution of
\begin{equation*}
a(z_h,\varphi_h) = \int_\Omega fz_h\dx + \int_\Gamma g z_h\dx\ \forall z_h\in Y_h.
\end{equation*}
Then
  \begin{align}
  \Vert \varphi-\varphi_h\Vert_{L^2(\Omega)} & + h^s \Vert \varphi-\varphi_h\Vert_{H^1(\Omega)}  \leq C h^{2s} \Vert \varphi\Vert_{W^{2,2}_{\vec\beta}(\Omega)}
  \nonumber\\
   &\leq C h^{2s} \left( \Vert f\Vert_{L^{2}_{\vec\beta}(\Omega)} + \sum_{j=1}^m \Vert g\Vert_{W^{1/2,2}_{\vec\beta}(\Gamma_j)}\right).
   \label{E4.15}
  \end{align}
\end{theorem}
\begin{proof}We will prove \eqref{E4.13} and \eqref{E4.14}. The proof of \eqref{E4.15} follows the same lines.

  We first prove that
  \begin{equation}\label{E4.16}\Vert y-y_h\Vert_{L^2(\Omega)}\leq C_{\mathcal{A}^*,\beta} \hat c_{\vec \mu} \Vert \mathcal{A}\Vert  h^s \Vert y-y_h\Vert_{H^1(\Omega)}\end{equation}
  Consider $\psi\in W^{2,2}_\beta(\Omega)$ the solution of the adjoint problem
  \[a(z,\psi) =\int_\Omega(y-y_h)z\mathrm{d}x\ \forall z\in H^1(\Omega)\]
  and let $\hat\psi_h\in Y_h$ be its Ritz-Galerkin projection onto $Y_h$ in the sense of $H^1(\Omega)$, as in the proof of Theorem \ref{T4.2}. We have, with \eqref{E4.6}, that
  \begin{align*}
    \Vert y-y_h\Vert_{L^2(\Omega)}^2 &= a(y-y_h,\psi) = a(y-y_h,\psi-\hat \psi_h) \\
    & \leq  \Vert \mathcal{A}\Vert  \Vert y-y_h\Vert_{H^1(\Omega)}\Vert \psi-\hat \psi_h\Vert_{H^1(\Omega)} \\
    &\leq  C_{\mathcal{A}^*,\beta} \hat c_{\vec \mu} \Vert \mathcal{A}\Vert  h^s \Vert y-y_h\Vert_{H^1(\Omega)}\Vert y-y_h\Vert_{L^2(\Omega)}
  \end{align*}
  and \eqref{E4.16} follows. Estimate \eqref{E4.13} follows from this, Theorem \ref{T4.2} and Lemma \ref{L2.6}.

  Using G\r{a}rding's inequality established in Lemma \ref{L2.3}, estimate \eqref{E4.16}, and the definition of $h_0>0$ in \eqref{E4.7}, we have that for all $h < h_0$
  \begin{align*}
    \frac{\Lambda}{8 C_E^2} \Vert y-y_h\Vert_{H^1(\Omega)}^2&\leq  a(y-y_h,y-y_h) + C_{\Lambda,E,b} \Vert y-y_h\Vert_{L^2(\Omega)}^2 \\
   &\leq  a(y-y_h,y-y_h) + C_{\Lambda,E,b} \Big(C_{\mathcal{A}^*,\beta} \hat c_{\vec \mu} \Vert \mathcal{A}\Vert  h^s\Big
    )^2 \Vert y-y_h\Vert_{H^1(\Omega)}^2\\
    &\leq a(y-y_h,y-y_h) + \frac{1}{4} \frac{\Lambda}{8 C_E^2} \Vert y-y_h\Vert_{H^1(\Omega)}^2,
  \end{align*}
  and hence
  \begin{equation}\label{E4.17}
    \frac{3\Lambda}{32 C_E^2} \Vert y-y_h\Vert_{H^1(\Omega)}^2\leq  a(y-y_h,y-y_h)
  \end{equation}
  Using Theorem \ref{T3.4}(c) and Lemma \ref{L4.1}
  \begin{equation}\label{E4.18}
    \Vert y- I_h y_h\Vert_{H^1(\Omega)}\leq \hat c_{\vec\mu} h^s \Vert y\Vert_{W^{2,2}_{\vec\beta}(\Omega)} \leq \hat c_{\vec\mu}C_{\mathcal{A},\vec\beta} h^s \sum_{j=1}^m\Vert u\Vert_{W^{1/2,2}_{\vec\beta}(\Gamma_j)}.
  \end{equation}
  Using that $a(y,I_h y_h) = a(y_h,I_h y_h)$, \eqref{E4.17} and the above inequality, we have that
  \begin{align*}
    \frac{3\Lambda}{32 C_E^2} \Vert y-y_h\Vert_{H^1(\Omega)}^2&\leq  a(y-y_h,y-I_h y_h)
    \leq   \Vert \mathcal{A}\Vert \Vert y-y_h\Vert_{H^1(\Omega)}\Vert y-I_h y_h\Vert_{H^1(\Omega)}\\
     &\leq \hat c_{\vec\mu}C_{\mathcal{A},\vec\beta} \Vert \mathcal{A}\Vert  h^s  \Vert y-y_h\Vert_{H^1(\Omega)},
  \end{align*}
  and the result follows.
\end{proof}

\begin{corollary}\label{C4.6}
There exists $C>0$ that depends on $A$, $b$, $a_0$, $\Omega$, and the mesh grading parameter $\vec\mu$, but is independent of $h$ such that
  for $0< h<\bar h$
  \begin{equation}\label{E4.19}
  \Vert y-y_h\Vert_{L^2(\Omega)} \leq C h^{3s/2} \Vert u\Vert_{L^{2}(\Gamma)}\ \forall u\in L^2(\Gamma)\end{equation}
   for all $s\leq 1$ and $s <\frac{\lambda_j}{\mu_j}$ for all $j\in \{1,\ldots,m\}$.

   Further, for all $f\in L^{2}_{\vec\beta}(\Omega)$ and $g\in \prod_{j=1}^m W^{1/2,2}_{\vec\beta}(\Gamma)$ and all $\theta\in(0,1)$, then
   \begin{equation}\label{E4.20}\Vert \varphi-\varphi_h\Vert_{H^{\theta}(\Omega)}\leq C h^{(2-\theta)s} \left( \Vert f\Vert_{L^{2}_{\vec\beta}(\Omega)} + \sum_{j=1}^m \Vert g\Vert_{W^{1/2,2}_{\vec\beta}(\Gamma_j)}\right),\end{equation}
   where $C$ is independent of $\theta$.
   \end{corollary}
\begin{proof}
  If $u\in H^{1/2}(\Gamma)$, then, by \eqref{E4.14} and the embedding $H^{1/2}(\Gamma)\hookrightarrow W^{1/2,2}_{\vec\beta}(\Gamma)\hookrightarrow\prod_{j=1}^m W^{1/2,2}_{\vec\beta}(\Gamma_j)$ for some $\vec\beta$ with $\beta_j\geq 0$, $1-\lambda_j < \beta_j < 1$, we obtain
  \[\Vert y-y_h\Vert_{L^2(\Omega)} \leq C h^{2s} \Vert u\Vert_{H^{1/2}(\Gamma)}.\]
  The first result follows by complex interpolation between this estimate and \eqref{E4.13}.

  The second one follows by interpolation between the estimates for $\theta = 0$ and $\theta=1$ that follow from \eqref{E4.15}.
\end{proof}
\section{\label{S5}Analysis of the control problem}

Now, we turn to the analysis of the control problem
\[
\Pb\quad \min_{u \in \uad} J(u) := \frac{1}{2}\int_\Omega (y_u(x) - y_d(x))^2\, \dx + \frac{\nu}{2}\int_\Gamma u^2(x)\, \dx +\int_\Gamma y_u(x) g_\varphi(x)\, \dx,
\]
where $y_u\in H^1(\Omega)$ solves \eqref{E2.7}.
For every $u\in H^{1/2}(\Gamma)'$, we define $\varphi_u\in H^1(\Omega)$ as the unique solution of
\[\langle z, \mathcal{A}^*\varphi_u\rangle_\Omega = \int_\Omega( y_u-y_d)z\dx +\int_\Gamma g_\varphi z\dx \forall z\in H^1(\Omega).
\]
We have that
\[J'(u)v = \int_\Omega(\varphi_u+\nu u)v\dx.\]
\begin{theorem}\label{T5.1}
  For any $y_d\in L^2(\Omega)$ and $g_\varphi \in L^2(\Gamma)$, problem \Pb  has a unique solution $\bar u\in\uad$ and there exist $\bar y,\bar\varphi\in H^1(\Omega)$ such that
  \begin{align*}
  \langle \mathcal{A}\bar y,z\rangle_\Omega &= \int_\Gamma  \bar u z\dx &&\forall z\in H^1(\Omega),\\
  \langle z, \mathcal{A}^*\bar \varphi\rangle_\Omega& = \int_\Omega(\bar y-y_d)z\dx +\int_\Gamma g_\varphi z\dx &&\forall z\in H^1(\Omega),\\
  \int_\Gamma(\bar\varphi+\nu\bar u)(u-\bar u)\dx&\geq 0 &&\forall u\in\uad,
\end{align*}
and $\bar u\in H^{1/2}(\Gamma)$.

If, further, $g_\varphi \in \prod_{j=1}^m W^{1/2,2}_{\vec\beta}(\Gamma_j)$ for some $\vec\beta$ such that $ 1-\lambda_j<\beta_j < 1$ and $\beta_j\ge0$ for all $j\in\{1,\ldots,m\}$,  then $\bar y, \bar\varphi\in W^{2,2}_{\vec\beta}(\Omega)\cap C(\bar\Omega)$, $\bar\varphi\in W^{3/2,2}_{\vec\beta}(\Gamma)\cap C(\Gamma)$, $\bar u\in C(\Gamma)$.

If, moreover, the weights also satisfy $\beta_j < 1/2$, for all $j\in\{1,\ldots,m\}$ then  $\bar\varphi, \bar u\in H^1(\Gamma)$.
\end{theorem}
\begin{proof}
  The existence of the solution follows from the appropriate continuity properties of the involved operators that are deduced from Lemma \ref{L2.6}. Uniqueness is deduced from the strict convexity of the functional. The first order optimality conditions are deduced, hence, in a standard way from the Euler-Lagrange equation $J'(\bar u)(u-\bar u)\geq 0$ for all $u\in\uad$ and Corollary \ref{C2.7}. The $H^1(\Omega)$ regularity of $\bar y$ follows from Lemma \ref{L2.3} and the regularity of the adjoint state from Lemma \ref{L2.6}. By the trace theorem, we have that $\varphi\in H^{1/2}(\Gamma)$. This and the projection formula
    \begin{equation}\label{E5.1}\bar u(x) = \proj_{[\umin,\umax]}\left(-\frac{\bar\varphi(x)}{\nu}\right),\end{equation}
  which follows in a standard way from the third optimality condition, imply the regularity of $\bar u$.

      Suppose now that $g_\varphi$ belongs to $L^2(\Gamma)\cap \prod_{j=1}^m W^{1/2,2}_{\vec\beta}(\Gamma_j)$ for some $\vec\beta$ such that $ 1-\lambda_j<\beta_j < 1$ and $\beta_j\ge0$ for all $j\in\{1,\ldots,m\}$.
  The $W^{2,2}_{\vec\beta}(\Omega)$ regularity of the state and adjoint state follow from a bootstrapping argument:
    since $\bar y\in H^1(\Omega)$ and $\beta_j\geq 0$ for all $j$, we have that $\bar y-y_d\in L^2(\Omega)\hookrightarrow  L^2_{\vec\beta}(\Omega)$.  From Theorem \ref{T3.5}(c) we deduce that $\bar\varphi\in W^{2,2}_{\vec\beta}(\Omega)$. This readily implies that $\bar\varphi\in W^{3/2,2}_{\vec\beta}(\Gamma)$. Using that $L^2_{\vec\beta}(\Omega)\subset L^r(\Omega)$ for all $1<r < 2/(1+\beta_j)$, we deduce from Theorem \ref{T3.5}(b) that $\bar\varphi\in W^{2,r}(\Omega)\hookrightarrow C(\bar\Omega)$, so $\bar\varphi\in C(\Gamma)$. Again the projection formula leads to $\bar u\in C(\Gamma)$.

  If $\beta_j < 1/2$, then $2/(1+\beta_j) > 4/3$, so there exists $r>4/3$ such that  $\bar\varphi\in W^{2,r}(\Omega)\hookrightarrow H^{3-2/r}(\Omega)$. Since $3-2/r>3/2$, by the trace theorem we have that $\bar\varphi\in C(\Gamma)\cap_{j=1}^m H^1(\Gamma_j) = H^1(\Gamma)$. This last equality follows because $\Gamma$ is one-dimensional and polygonal. This regularity is preserved by the projection formula, and therefore $\bar u\in H^1(\Gamma)$.
  \end{proof}

Notice that for any polygonal domain $\lambda_j > 1/2$ for all $j\in\{1,\ldots,m\}$, so the condition $\beta_j < 1/2$ may be a constraint in the regularity of the datum $g_\varphi$, but it is not a constraint on the domain.
Although some of the intermediate results below can be proved for $g_\varphi\in L^2(\Gamma)$, since the main result requires $H^1(\Gamma)$ regularity of the optimal control, in the rest of the work we will do the following assumption.
\begin{assumption}\label{A5.2}We assume that $g_\varphi \in \prod_{j=1}^m W^{1/2,2}_{\vec\beta}(\Gamma_j)$ for some $\vec\beta$ such that $1-\lambda_j<\beta_j<1/2$, $\beta_j\geq 0$ for all $j\in\{1,\ldots,m\}$.  We denote
\[M_d = \Vert y_d\Vert_{L^2(\Omega)} + \sum_{j =1}^m\Vert g_\varphi\Vert_{W^{1/2,2}_{\vec\beta}(\Gamma_j)} + 1.\]
\end{assumption}

For every $u\in L^2(\Gamma)$, we will denote $y_h(u)$ the solution of the discrete state equation \eqref{E4.1} and $\varphi_h(u)$ the solution of
\[a(z_h,\varphi_h) = \int_\Omega (y_h(u)-y_d) z_h\dx \int_\Gamma g_\varphi z_h\dx\ \forall z_h\in Y_h.\]
Our discrete functional reads like
\[J_h(u) = \frac{1}{2}\int_\Omega (y_h(u)-y_d)^2\dx + \frac{\nu}{2}\int_\Gamma u^2\dx+\int_\Gamma y_h(u) g_\varphi\dx.\]
To discretize the control, we notice that every triangulation $\mathcal{T}_h$ of $\Omega $ defines a segmentation $\mathcal{E}_h$ of $\Gamma$ and define $\uadh = \Uh\cap \uad$, where
\[\Uh = \{u_h\in L^2(\Gamma):\ u_{h\vert E}\in \mathcal{P}^0(E)\ \forall E\in\mathcal{E}_h\}.\]
Here and elsewhere $\mathcal{P}^i(K)$ is the set of polynomials of degree $i$ in the set $K$.
For every $u\in L^1(\Gamma)$, we define $Q_h u\in\Uh$ by
\[Q_h u(x) = \displaystyle\frac{1}{h_E}  \int_E u\dx  \mbox{ if }x\in E,\]
where $E\in\mathcal{E}_h$ and $h_E$ is the length of $E$. Notice that $u\in\uad$ implies $Q_h u\in\uadh$.
\begin{lemma}\label{L5.3}
For every $u\in H^1(\Gamma)$ there exists a constant $C>0$ independent of $h$ such that
  \[\Vert u-Q_h u\Vert_{(H^{1}(\Gamma))'}+ h \Vert  u-Q_h u\Vert_{L^2(\Gamma)} \leq C h^2\Vert u\Vert_{H^1(\Gamma)}.\]
If Assumption \ref{A5.2} holds, then we also have that
\[\left\vert\int_\Gamma(\varphi_u + \nu u)(u-Q_h u)\dx\right\vert\leq C h^2\left(\Vert u\Vert_{H^1(\Gamma)}^2+ M_d^2 \right).\]
\end{lemma}
\begin{proof}
It is well known that for every $E\in\mathcal{E}_h$ we have $\Vert  u-Q_h u\Vert_{L^2(E)} \leq C h_E\Vert u\Vert_{H^1(E)}$. Using that $h_E\leq c h$, we have
\[\Vert  u-Q_h u\Vert_{L^2(\Gamma)}^2 = \sum_{E\in\mathcal{E}_h} \Vert  u-Q_h u\Vert_{L^2(E)}^2 \leq C \sum_{E\in\mathcal{E}_h}h_E^2\Vert u\Vert_{H^1(E)}^2 \leq C h^2 \Vert u\Vert_{H^1(\Gamma)}^2.\]
The estimate for the norm in $H^1(\Gamma)'$ follows now by duality since $\int_\Gamma(u-Q_hu)w_h\dx=0$ for all $w_h\in\Uh$. This estimate implies the third one taking into account that, using the same arguments as in the proof of Theorem \ref{T5.1},  $\varphi_u\in H^1(\Gamma)$, and
\begin{align*}
  \Vert &\varphi_u\Vert_{H^1(\Gamma)}  \leq  C \Vert \varphi\Vert_{W^{2,2}_{\vec\beta}(\Omega)}
  \leq  C \left(\Vert y_u-y_d\Vert_{L^2_{\vec\beta}(\Omega)} + \sum_{j = 1}^m \Vert g_\varphi\Vert_{W^{1/2,2}_{\vec\beta}(\Gamma_j)}\right) \\
  &\leq  C \left(\Vert y_u\Vert_{L^2(\Omega)}+ \Vert y_d\Vert_{L^2(\Omega)} + \sum_{j = 1}^m \Vert g_\varphi\Vert_{W^{1/2,2}_{\vec\beta}(\Gamma_j)}\right)   \leq  C \left(\Vert u\Vert_{L^2(\Gamma)}+ M_d\right).
\end{align*}
Therefore, we obtain
\begin{align*}
  \left\vert\int_\Gamma(\varphi_u + \nu u)(u-Q_h u)\dx\right\vert & \leq  \Vert \varphi_u + \nu u\Vert_{H^1(\Gamma)} \Vert u-Q_h u\Vert_{H^1(\Gamma)'}\\
  &\leq C \left(M_d + \Vert u\Vert_{L^2(\Gamma)} + \nu \Vert u\Vert_{H^1(\Gamma)}\right) h^2 \Vert u\Vert_{H^1(\Gamma)}
\end{align*}
and the result follows using Young's inequality.
\end{proof}
Our discrete problems reads like
\[\Pbh\ \min_{u_h\in\uadh} J_h(u_h).\]
Existence and uniqueness of solution of problem \Pbh, as well as first order optimality conditions follow in an standard way. We state the result in the next theorem for further reference.
\begin{theorem}
  For every $0<h<\bar h$, problem \Pbh has a unique solution $\bar u_h\in\uadh$. Further, if we denote $\bar y_h = y_h(\bar u_h)$ and $\bar\varphi_h = \varphi_h(\bar u_h)$, then
  \begin{equation}\label{E5.2}\int_\Gamma(\bar\varphi_h + \nu\bar u_h)(u_h-\bar u_h)\dx \geq 0\ \forall u_h\in\uadh.\end{equation}
\end{theorem}

Before stating and proving the main theorem of this section, we prove two auxiliary results.

\begin{lemma}\label{L5.5}There exists $C>0$ independent of $h$, $y_d$ and $g_\varphi$ such that for all $0<h<\bar h$,
\[\Vert \bar y\Vert_{H^1(\Omega)}+\Vert \bar u\Vert_{H^{1/2}(\Gamma)} +
\Vert \bar y_h\Vert_{H^1(\Omega)}+\Vert \bar u_h\Vert_{L^2(\Gamma)} \leq C \left(\Vert y_d\Vert_{L^2(\Omega)} + \Vert g_\varphi\Vert_{L^2(\Gamma)}+1\right).
\]
If, moreover, Assumption \ref{A5.2} holds, then
\[\Vert \bar u\Vert_{H^1(\Gamma)}\leq C M_d.\]
\end{lemma}

\begin{proof}Consider a fixed  $u_{\mathrm{ad}}\in\uad$ such that $u_{\mathrm{ad}}\in\uadh$ for all $h>0$.
    Using that $\Vert \bar y_h-y_d\Vert^2_{L^2(\Omega)}\geq 0$ and the optimality of $\bar u_h$ together with Young's inequality and estimate \eqref{E4.12}, we have for all $\varepsilon > 0$ that
  \begin{align*}
  \frac{\nu}{2}&\Vert \bar u_h\Vert^2_{L^2(\Gamma)}\leq  J_h(\bar u_h)- \int_\Gamma y_h(\bar u_h)g_\varphi\dx \\
   & \leq J_h(u_{\mathrm{ad}}) + \varepsilon \Vert  y_h(\bar u_h)\Vert^2_{L^2(\Gamma)} + \frac{1}{4\varepsilon} \Vert g_\varphi\Vert_{L^2(\Gamma)}^2\\
  &\leq  \frac{1}{2}\Vert y_h(u_{\mathrm{ad}}) - y_d\Vert_{L^2(\Omega)}^2 + \frac{\nu}{2}\Vert u_{\mathrm{ad}}\Vert^2_{L^2(\Gamma)} + \int_\Gamma y_h(u_{\mathrm{ad}}) g_\varphi\dx \\
  & \qquad + \varepsilon c_2^2 \Vert \bar u_h\Vert^2_{L^2(\Gamma)} +  \frac{1}{4\varepsilon} \Vert g_\varphi\Vert_{L^2(\Gamma)}^2\\
  & \leq \Vert y_h(u_{\mathrm{ad}})\Vert_{L^2(\Omega)}^2+ \Vert  y_d\Vert_{L^2(\Omega)}^2 + (\frac{\nu}{2}+c_2^2)\Vert u_{\mathrm{ad}}\Vert^2_{L^2(\Gamma)} \\
  & \qquad
  + \varepsilon c_2^2 \Vert \bar u_h\Vert^2_{L^2(\Gamma)} +  \frac{1+\varepsilon}{4\varepsilon} \Vert g_\varphi\Vert_{L^2(\Gamma)}^2
  \end{align*}
where $c_2$ is introduced in \eqref{E4.12}.
  Taking $\varepsilon = \nu/(4c_2^2)$, we readily deduce that  $\{\bar u_h\}$ is uniformly bounded in $L^2(\Gamma)$.
   The estimate for $\Vert \bar y_h\Vert_{H^1(\Omega)}$ follows from this one and estimate \eqref{E4.2}.

  Estimates for $\Vert \bar u\Vert_{L^2(\Gamma)}$ and $\Vert \bar y\Vert_{H^1(\Omega)}$ follow in a similar way. From this last one and Lemma \ref{L2.6} an estimate for $\Vert \bar\varphi\Vert_{H^1(\Omega)}$ in terms of the data is obtained. The trace theorem and the projection formula \eqref{E5.1} lead to the estimate for $\Vert \bar u\Vert_{H^{1/2}(\Gamma)}$.

  If Assumption \ref{A5.2} holds, then, using the estimate for $\Vert \bar y\Vert_{L^2(\Omega)}$ and noting that the condition $\beta_j<1/2$ implies $\prod_{j=1}^m W^{1/2,2}_{\vec\beta}(\Gamma_j)\hookrightarrow L^2(\Gamma)$ and hence
\[\Vert y_d\Vert_{L^2(\Omega)} + \Vert g_\varphi\Vert_{L^{2}(\Gamma)} + 1\leq M_d,\]
 we obtain an estimate of $\Vert \bar\varphi\Vert_{W^{2,2}_{\vec\beta}(\Omega)}$ in terms of $M_d$. The trace theorem and the projection formula \eqref{E5.1} lead to the estimate for $\Vert \bar u\Vert_{H^1(\Gamma)}$.
\end{proof}
In the rest of the work $s$ represents any positive number satisfying $s\leq 1$ and $s<\lambda_j/\mu_j$.

\begin{lemma}\label{L5.6}Suppose Assumption \ref{A5.2} holds. Then, there exists $C>0$ independent of $h$, $y_d$, $g_\varphi$ and $\{\bar u_h\}$ such that
  \begin{equation}\label{E5.3}
    \Vert \varphi_{\bar u_h}-\bar\varphi_h\Vert_{L^2(\Omega)}\leq  C h^{3s/2} M_d.
  \end{equation}
Moreover, for all $\theta\in (1/2,1]$ we have the following estimate:
\begin{equation}\label{E5.4}  \Vert  \varphi_{\bar u_h}-\bar \varphi_h \Vert_{L^2(\Gamma)} \leq C h^{(2-\theta)s}M_d.  \end{equation}
\end{lemma}
\begin{proof}
By the triangle inequality
\begin{equation}\label{E5.5}\Vert \varphi_{\bar u_h}-\bar\varphi_h\Vert_{L^2(\Omega)}\leq \Vert \varphi_{\bar u_h}-\varphi^h\Vert_{L^2(\Omega)} + \Vert \varphi^h-\bar\varphi_h\Vert_{L^2(\Omega)},\end{equation}
where $\varphi^h$ is the unique element in $H^1(\Omega)$ such that $a(z,\varphi^h) =  \int_\Omega(\bar y_h-y_d)z\dx +  \int_\Gamma g_\varphi z\dx$ for all $z\in H^1(\Omega)$, i.e., $\bar\varphi_h$ is the finite element approximation of $\varphi^h$.

Let us estimate the first term in the right hand side of \eqref{E5.5}.
Noting that
  \[a(z,\varphi_{\bar u_h}-\varphi^h) = \int_\Omega(y_{\bar u_h}-y_h(\bar u_h))(z)\dx\ \forall z\in H^1(\Omega),\]
  we deduce from Theorem \ref{T3.5},  the existence of $C>0$ independent of $h$ such that
  \begin{equation}\label{E5.6}\Vert \varphi_{\bar u_h}-\varphi^h\Vert_{L^2(\Omega)}\leq C \Vert y_{\bar u_h}-y_h(\bar u_h)\Vert_{L^2(\Omega)}.\end{equation}
Applying the finite element error estimate for the state \eqref{E4.19} of Corollary \ref{C4.6} and Lemma \ref{L5.5}, we have
  \begin{align*}
  \Vert y_{\bar u_h}-y_h(\bar u_h)\Vert_{L^2(\Omega)} & \leq C h^{3s/2}\Vert \bar u_h\Vert_{L^2(\Gamma)} \\
  & \leq C h^{3s/2}(\Vert y_d\Vert_{L^2(\Omega)} + \Vert g_\varphi\Vert_{L^2(\Gamma)}+1)\leq C h^{3s/2} M_d.
  \end{align*}
This, together with \eqref{E5.6} leads to
  \begin{equation}\label{E5.7}
  \Vert \varphi_{\bar u_h}-\varphi^h\Vert_{L^2(\Omega)}\leq C h^{3s/2}M_d.
  \end{equation}
  To estimate the second summand in the right hand side of \eqref{E5.5} we apply the finite element error estimate \eqref{E4.15}, the uniform boundness result in Lemma \ref{L5.5} and the embedding $\prod_{j = 1}^m W^{1/2,2}_{\vec\beta}(\Gamma) \hookrightarrow L^2(\Gamma)$:
  \begin{align*}
    \Vert \varphi^h-\bar\varphi_h\Vert_{L^2(\Omega)} &\leq  C \left(\Vert \bar y_h-y_d\Vert_{L^2(\Omega)} + \sum_{j =1}^m\Vert g_\varphi\Vert_{W^{1/2,2}_{\vec\beta}(\Gamma_j)}\right) h^{2s}\\
    &\leq C \left(\Vert \bar y_h\Vert_{L^2(\Omega)}+ \Vert y_d\Vert_{L^2(\Omega)} + \sum_{j =1}^m\Vert g_\varphi\Vert_{W^{1/2,2}_{\vec\beta}(\Gamma_j)}  \right) h^{2s}\\
    &\leq C \left(2\Vert y_d\Vert_{L^2(\Omega)} + \Vert g_\varphi\Vert_{L^2(\Gamma)} + \sum_{j =1}^m\Vert g_\varphi\Vert_{W^{1/2,2}_{\vec\beta}(\Gamma_j)} +1 \right) h^{2s}\\
    &\leq C \left(\Vert y_d\Vert_{L^2(\Omega)} + \sum_{j =1}^m\Vert g_\varphi\Vert_{W^{1/2,2}_{\vec\beta}(\Gamma_j)} +1 \right) h^{2s} = C h^{2s}M_d.
  \end{align*}
  Estimate \eqref{E5.3} follows, hence, from \eqref{E5.5} together with this last estimate and \eqref{E5.7}.

Let us prove \eqref{E5.4}. First we notice that for $1/2<\theta\leq 1$, the trace operator is continuous from $H^\theta(\Omega)$ to $L^2(\Gamma)$, so
  \begin{align*}
    \Vert  \varphi_{\bar u_h}-\bar \varphi_h \Vert_{L^2(\Gamma)} &\leq C \Vert  \varphi_{\bar u_h}-\bar \varphi_h \Vert_{H^\theta(\Omega)} .
  \end{align*}
 To estimate the term $\Vert  \varphi_{\bar u_h}-\bar \varphi_h \Vert_{H^\theta(\Omega)}$, we first introduce $\phi_h\in Y_h$, the finite element approximation of $\varphi_{\bar u_h}$, that satisfies $a(z_h,\phi_h) =   \int_\Omega( y_{\bar u_h}-y_d)z_h\dx +  \int_\Gamma g_\varphi z_h\dx$ for all $z_h\in Y_h$.
The difference $\phi_h-\bar\varphi_h$ satisfies $a(z_h,\phi_h-\bar\varphi_h) = \int_\Omega (y_{\bar u_h}-\bar y_h) z_h\dx$ for all $z_h\in Y_h$. From the continuity estimate for the discrete adjoint equation of Theorem \ref{T4.3} we deduce that
\begin{equation}\label{E5.8}\Vert  \phi_h-\bar \varphi_h \Vert_{H^1(\Omega)}\leq C\Vert y_{\bar u_h}- \bar y_h\Vert_{L^2(\Omega)}.  \end{equation}
   Using the triangle inequality, the fact that $\theta\leq 1$, the finite element error estimate for the adjoint estate equation  \eqref{E4.19} of Corollary \ref{C4.6}, \eqref{E5.8}, and the finite element error estimate for the state equation \eqref{E4.20}, together with the uniform boundness of $\Vert \bar u_h\Vert_{L^2(\Gamma)}$ provided in Lemma \ref{L5.5}, we obtain
  \begin{align*}
    \Vert  \varphi_{\bar u_h}-\bar \varphi_h \Vert_{H^\theta(\Omega)}
    &\leq   \Vert  \varphi_{\bar u_h}-\phi_h \Vert_{H^\theta(\Omega)} + \Vert  \phi_h-\bar \varphi_h \Vert_{H^1(\Omega)}\\
    &\leq     C\left(h^{(2-\theta)s} M_d +  \Vert y_{\bar u_h}- \bar y_h\Vert_{L^2(\Omega)}     \right)\\
    &\leq C\left(h^{(2-\theta)s} M_d +  h^{3s/2}\Vert \bar u_h\Vert_{L^2(\Gamma)}     \right) \leq C h^{(2-\theta)s} M_d,
  \end{align*}
where the last inequality is a result of Lemma \ref{L5.5} and the condition $\theta > 1/2$.
\end{proof}
As already mentioned in the introduction, a straightforward application of the usual techniques of proof for the following theorem would lead only to order $h^s$ with $s\leq 1$ and $s<\lambda_j/\mu_j$. The next result shows an improvement of the order of convergence for the control variable with respect to the known results obtained for problems governed by coercive equations; cf. \cite{APR2012,APR2015}.
\begin{theorem}\label{T5.7}Suppose Assumption \ref{A5.2} holds. Then, there exists a constant independent of $h$, $y_d$ and $g_\varphi$ such that, for all $0< h<\bar h$
  \[\Vert \bar u-\bar u_h\Vert_{L^2(\Gamma)} \leq C h^{s^*} M_d,\]
for all $s^*\leq 1$ such that $s^* <\dfrac{3}{2}\dfrac{\lambda_j}{\mu_j}$ for all $j\in \{1,\ldots,m\}$.
\end{theorem}

\begin{proof}
  Testing the equality $a(z,\bar\varphi-\varphi_{\bar u_h}) =  \int_\Omega(\bar y-y_{\bar u_h})z\dx$ for $z = \bar y-y_{\bar u_h}$ and using the state equation, we have that
  \[0\leq \Vert \bar y-y_{\bar u_h}\Vert^2_{L^2(\Omega)} = a(\bar y-y_{\bar u_h},\bar\varphi-\varphi_{\bar u_h}) = \int_\Gamma(\bar u-\bar u_h)(\bar\varphi-\varphi_{\bar u_h})\dx.\]
  So we can write
  \begin{align*}
    \nu\Vert \bar u&-\bar u_h\Vert_{L^2(\Gamma)}^2 \leq  \int_\Gamma(\bar\varphi-\varphi_{\bar u_h} + \nu (\bar u-\bar u_h))( \bar u-\bar u_h)\dx \\
     &=  \int_\Gamma(\bar\varphi-\bar \varphi_h + \nu (\bar u-\bar u_h))(\bar u-\bar u_h)\dx + \int_\Gamma(\bar\varphi_h-\varphi_{\bar u_h} )( \bar u-\bar u_h)\dx = I+II.
  \end{align*}
   Let us bound the first term.
    First we insert in appropriate places $Q_h\bar u$ and $\bar u$. Next, we apply the first order optimality conditions for the continuous and discrete problem. Finally we insert $\varphi_{\bar u_h}$ to obtain
  \begin{align*}
    I &=   \int_\Gamma(\bar\varphi-\bar \varphi_h + \nu (\bar u-\bar u_h) )( \bar u- Q_h\bar u)\dx +
    \int_\Gamma(\bar\varphi-\bar \varphi_h + \nu (\bar u-\bar u_h) )( Q_h\bar u-\bar u_h)\dx \\
    &=   \int_\Gamma(\bar\varphi-\bar \varphi_h + \nu (\bar u-\bar u_h) )( \bar u- Q_h\bar u)\dx +
    \int_\Gamma(\bar\varphi + \nu \bar u )( Q_h\bar u-\bar u_h)\dx \\
    &\quad +
    \int_\Gamma(\bar\varphi_h + \nu \bar u_h) )( \bar u_h-Q_h\bar u)\dx \\
    &=   \int_\Gamma(\bar\varphi-\bar \varphi_h + \nu (\bar u-\bar u_h) )( \bar u- Q_h\bar u)\dx +
    \int_\Gamma(\bar\varphi + \nu \bar u )( Q_h\bar u-\bar u)\dx \\
    & \quad+
    \int_\Gamma(\bar\varphi + \nu \bar u )( \bar u-\bar u_h)\dx +
    \int_\Gamma(\bar\varphi_h + \nu \bar u_h) )( \bar u_h-Q_h\bar u)\dx \\
   &\leq  \int_\Gamma(\bar\varphi-\bar \varphi_h + \nu (\bar u-\bar u_h) )( \bar u- Q_h\bar u)\dx
       +\int_\Gamma (\bar\varphi + \nu \bar u )( Q_h\bar u-\bar u)\dx \\
   &=  \int_\Gamma(\bar\varphi-\varphi_{\bar u_h} + \nu (\bar u-\bar u_h) )( \bar u- Q_h\bar u)\dx
   + \int_\Gamma(\varphi_{\bar u_h}-\bar \varphi_h )( \bar u- Q_h\bar u)\dx\\
   &\quad +    \int_\Gamma(\bar\varphi + \nu \bar u )( Q_h\bar u-\bar u)\dx  =  I_A + I_B + I_C.
  \end{align*}
  From Lemmas \ref{L5.3} and \ref{L5.5}, it is clear that $I_C\leq C h^2 M_d^2$.

Let us study $I_A$.  Testing the equality $a(z,\bar\varphi-\varphi_{\bar u_h}) = \int_\Omega(\bar y-y_{\bar u_h})z\dx$ for $z = \bar y-y_{Q_h\bar u}$ and using the state equation, Cauchy-Schwarz inequality, and Theorem \ref{T3.4}(a),  we obtain
  \begin{align*}
  \int_\Gamma(\bar\varphi & -\varphi_{\bar u_h})(\bar u- Q_h\bar u)\dx =
   a(\bar y-y_{Q_h\bar u},\bar\varphi-\varphi_{\bar u_h})
   =  \int_\Omega(\bar y-y_{\bar u_h})(\bar y-y_{Q_h\bar u})\dx\\
  &\leq  \Vert y_{\bar u-Q_h\bar u}\Vert_{L^2(\Omega)} \Vert y_{\bar u-\bar u_h}\Vert_{L^2(\Omega)} \leq C \Vert \bar u-Q_h\bar u\Vert_{L^2(\Gamma)} \Vert \bar u-\bar u_h\Vert_{L^2(\Gamma)}
  \end{align*}
  Using this and Lemmas \ref{L5.3} and \ref{L5.5}, we obtain
  \[
  I_A = \int_\Gamma(\bar\varphi-\varphi_{\bar u_h} + \nu (\bar u-\bar u_h) )( \bar u- Q_h\bar u)\dx\leq C h \Vert \bar u-\bar u_h\Vert_{L^2(\Gamma)}M_d.
   \]
   Next we bound $I_B$ and $II$. By the Cauchy-Schwarz inequality, we have that, for every $v\in L^2(\Gamma)$,
   \begin{equation}\label{E5.9}
   \int_\Gamma(\varphi_{\bar u_h}-\bar \varphi_h)v\dx\leq   \Vert  \varphi_{\bar u_h}-\bar \varphi_h\Vert_{L^2(\Gamma)} \Vert v\Vert_{L^2(\Gamma)}.
   \end{equation}
  Taking $v =  \bar u- Q_h\bar u$ in \eqref{E5.9} and using \eqref{E5.4} and Lemmas \ref{L5.3} and \ref{L5.5}, we  conclude that

  \begin{align*} I_B \leq  \Vert  \varphi_{\bar u_h}-\bar \varphi_h\Vert_{L^2(\Gamma)} \Vert  \bar u- Q_h\bar u\Vert_{L^2(\Gamma)} \leq C h^{(2-\theta)s+1} M_d^2.
  \end{align*}

  Finally, taking $v = \bar u-\bar u_h$ in \eqref{E5.9} and using \eqref{E5.4}, we have

\begin{align*}
    II \leq
      C h^{(2-\theta)s} \Vert \bar u-\bar u_h\Vert_{L^2(\Gamma)} M_d.
  \end{align*}

  Gathering all the estimates we have that

  \begin{align*}
  \nu \Vert \bar u-\bar u_h\Vert_{L^2(\Gamma)}^2 &\leq C( h \Vert \bar u-\bar u_h\Vert_{L^2(\Gamma)} M_d + h^{(2-\theta)s+1} M_d^2 \\
  & \qquad + h^2 M_d^2 + h^{(2-\theta)s}\Vert \bar u-\bar u_h\Vert_{L^2(\Gamma)} M_d)
  \end{align*}
  and the proof concludes  using Young's inequality.
   Notice that the appearance of the terms $h \Vert \bar u-\bar u_h\Vert_{L^2(\Gamma)} M_d $ and  $h^2 M_d^2$ implies that the resulting exponent $s^*$ is less or equal than one.
   On the other hand, since $\theta > 1/2$, the term  $h^{(2-\theta)s}\Vert \bar u-\bar u_h\Vert_{L^2(\Gamma)} M_d$ yields the bound $s^*\leq (2-\theta)s < \dfrac{3}{2}s <\dfrac{3}{2}\dfrac{\lambda_j}{\mu_j}$. Finally, from the term $h^{(2-\theta)s+1} M_d^2$  we obtain the bound $s^*\leq \min\{(2-\theta)s,1\}$, so no new conditions are imposed on $s^*$.
\end{proof}

\section{A numerical example}\label{S6}
Let $\Omega$ be the $L$-shaped domain  $\Omega = \{x\in\mathbb{R}^2: r<\sqrt{2}, \theta < 3\pi/2\}\cap (-1,1)^2$.
We consider a functional of the form
\[J(u) = \frac{1}{2}\int_\Omega (y_u(x)-y_d(x))^2 \dx +\frac{\nu}{2}\int_\Gamma u(x)^2 \dx+\int_\Gamma y_u(x) g_\varphi(x) \dx,\]
where
\begin{equation*}
\left\{\begin{array}{rcl} -\Delta y_u + b\cdot\nabla y_u + a_0 y_u &=& f \text{ in } \Omega,\\ \partial_{n} y &=& u+g_y\text{ on } \Gamma.\end{array}\right.
\end{equation*}
with data $\nu$, $y_d$, $g_\varphi$, $b$, $a_0$, $g_y$ described below. The inclusion of data $f$ and $g_y$ is useful to write a problem with known exact solution. Notice that, if we denote $y_0\in L^2(\Omega)$ the state related to $u\equiv 0$ and redefine $y_d:=y_d-y_0$ and $y_u:=y_u-y_0$, the problem fits into the framework of problem \Pb and equation \eqref{E1.1}.

Let $(r,\theta)$ be the polar coordinates in the plane, $r\geq 0$, $\theta\in [0,2\pi]$. The interior angle at the vertex of the domain located at the origin is $\omega = \omega_1 = 3\pi/2$ and we denote $\lambda = \lambda_1 =\pi/\omega_1 = 2/3$. For $j=2,\ldots,6$, $\omega_j=\pi/2$ and $\lambda_j=2$.

We introduce $\bar y = r^\lambda\cos(\lambda\theta),$ $\bar\varphi = -\bar y$ and $\bar u = -\bar\varphi/\nu$ on $\Gamma$ and, for some $\alpha> -3/2$ and some $\delta \geq 0$, we consider $b(x) = \delta r^{\alpha+1}(\cos\theta,\sin\theta)^T$ and $a_0(x) = r^\alpha$.

The data for this problem are defined as $f = b\cdot\nabla \bar y+ a_0 \bar y$, $ g_y = \partial_{n} \bar y - \bar u$ on $\Gamma$,  $ y_d = \bar y +\nabla\cdot(\bar\varphi b)-a_0\bar\varphi$ and $g_\varphi = \partial_n\bar\varphi + (b\bar\varphi)\cdot n$.

For all $\alpha > -2$, $b\in L^\ps(\Omega)$ for some $\ps>2$ (Assumption \ref{A2.1}). For $\alpha > -1-\beta$, $a_0,\nabla\cdot b, f, y_d\in L^2_{\vec\beta}(\Omega)$ and $b\cdot n, g_y, g_\varphi \in W^{1/2,2}_{\vec\beta}(\Gamma)$, so the assumptions of theorems \ref{T3.4}(c) and \ref{T3.5}(c) hold. If we impose $\beta<1/2$ (assumption in Theorem \ref{T5.1}), we have that for $\alpha > -3/2$ all the assumptions of the paper hold. In our experiments, we fix $\alpha = -1.25$.

The given $\bar u$ is the solution of the control problem
\[\Pb \min_{u\in L^2(\Gamma)}J(u),\]
with related state $\bar y$ and adjoint state $\bar\varphi$, which satisfy the optimality system
\begin{equation*}
\left\{\begin{array}{rcl} -\Delta\bar y + b\cdot\nabla \bar y + a_0 \bar y &=& f \text{ in } \Omega,\\ \partial_{n_A} \bar y &=& g_y+\bar u\text{ on } \Gamma,\end{array}\right.
\end{equation*}
\begin{equation*}
\left\{\begin{array}{rccl} -\Delta\bar \varphi -\nabla\cdot(b \bar \varphi) + a_0 \bar \varphi &=& \bar y-y_d & \text{ in } \Omega,\\ \partial_{n} \bar \varphi + \bar \varphi b\cdot n &=& g_\varphi & \text{ on } \Gamma,\end{array}\right.
\end{equation*}
\begin{equation*}
  \bar u = -\bar\varphi/\nu\mbox{ on }\Gamma.
\end{equation*}
It is clear that  $\bar y,\bar\varphi\in W^{2,2}_{\vec\beta}(\Omega)$ and $\bar u\in H^1(\Gamma)\cap W^{1/2,2}_{\vec\beta}(\Gamma)$ for  $\vec\beta=(\beta,0,0,0,0,0)$ for all   $\beta > 1-\lambda >1/3$.

For $\delta = 6$, we have checked numerically that the operator is not coercive,

To discretize the problem we use the finite element approximation described in the work. We use a family of graded meshes obtained by bisection; see, e.g., \cite[Figure 1.2]{AMPR2019}. This meshing method does not lead to superconvergence properties in the gradients.

 First we check estimates \eqref{E4.14} and \eqref{E4.15} for the error in the solution of the boundary value problem. For appropriately graded meshes, $\mu < 2/3 = \lambda$, we expect order $h^2$ in $L^2(\Omega)$ and order $h$ in $H^1(\Omega)$. For a quasi-uniform family, $\mu = 1$, we have $s < 2/3$, so we expect order $h^{1.33}$ in $L^2(\Omega)$ and order $h^{0.66}$ in $H^1(\Omega)$. We summarize the results in tables \ref{T6.1}--\ref{T6.4}. We include results for both the state and adjoint state equation. Notice that $\tilde\varphi_h$ is the finite element approximation of $\bar\varphi$, obtained using the exact $\bar y$, i.e., $a(z_h,\tilde\varphi_h) =   \int_\Omega(\bar y-y_d)z_h\dx +  \int_\Gamma g_\varphi z_h\dx$ for all $z_h\in Y_h$.

\begin{table}
  \centering
  \[
  \begin{array}{c@{\;}|@{\;}c@{\;}c@{\;}|@{\;}c@{\;}c@{\;}||@{\;}c@{\;}c@{\;}|@{\;}c@{\;}c}
 j &  \Vert \bar y-\bar y_h(\bar u)\Vert_{L^2(\Omega)} &   EOC  &     \Vert \bar y-\bar y_h(\bar u)\Vert_{H^1(\Omega)}  &  EOC  &    \Vert \bar \varphi-\tilde \varphi_h\Vert_{L^2(\Omega)} &   EOC  &     \Vert \bar \varphi-\tilde \varphi_h\Vert_{H^1(\Omega)}  &  EOC  \\ \hline
 1 & 1.20e-01 &         &  2.92e-01 &           &  3.85e-02 &         &  3.28e-01 &        \\
 2 & 5.67e-02 &   1.08  &  1.94e-01 &   0.59    &  1.27e-02 &   1.60  &  2.02e-01 &   0.70 \\
 3 & 2.57e-02 &   1.14  &  1.25e-01 &   0.63    &  4.45e-03 &   1.51  &  1.27e-01 &   0.68 \\
 4 & 1.12e-02 &   1.20  &  7.98e-02 &   0.65    &  1.64e-03 &   1.44  &  7.98e-02 &   0.67 \\
 5 & 4.72e-03 &   1.25  &  5.05e-02 &   0.66    &  6.24e-04 &   1.39  &  5.04e-02 &   0.66 \\
 6 & 1.94e-03 &   1.28  &  3.19e-02 &   0.66    &  2.42e-04 &   1.37  &  3.19e-02 &   0.66 \\
 7 & 7.88e-04 &   1.30  &  2.02e-02 &   0.66    &  9.47e-05 &   1.35  &  2.01e-02 &   0.66 \\
 8 & 3.17e-04 &   1.32  &  1.27e-02 &   0.66    &  3.73e-05 &   1.34  &  1.27e-02 &   0.66 \\
 9 & 1.27e-04 &   1.32  &  8.02e-03 &   0.67    &  1.47e-05 &   1.34  &  8.01e-03 &   0.66 \\ \hline
\multicolumn{2}{l}{\text{Expected}}  &  1.33  &            &  0.66       &   &  1.33  &            &  0.66
 \end{array}
 \]
  \caption{Errors and experimental orders of convergence for the boundary value problem. In the bottom line, orders of convergence expected from \eqref{E4.14}. Quasi-uniform mesh family. }\label{T6.1}
\end{table}

\begin{table}
  \centering
 \[
  \begin{array}{c@{\;}|@{\;}c@{\;}c@{\;}|@{\;}c@{\;}c@{\;}||@{\;}c@{\;}c@{\;}|@{\;}c@{\;}c}
 j &    \Vert \bar y-\bar y_h(\bar u)\Vert_{L^2(\Omega)} &   EOC  &     \Vert \bar y-\bar y_h(\bar u)\Vert_{H^1(\Omega)}  &  EOC  &    \Vert \bar \varphi-\tilde \varphi_h\Vert_{L^2(\Omega)} &   EOC  &     \Vert \bar \varphi-\tilde \varphi_h\Vert_{H^1(\Omega)}  &  EOC \\ \hline
 1 &     8.30e-02 &          &  2.48e-01 &         &  2.27e-02 &         &  2.65e-01 &         \\
 2 &     2.59e-02 &   1.68   &  1.40e-01 &   0.82  &  5.53e-03 &   2.03  &  1.42e-01 &   0.90  \\
 3 &     8.07e-03 &   1.68   &  7.74e-02 &   0.86  &  1.50e-03 &   1.89  &  7.76e-02 &   0.87  \\
 4 &     2.49e-03 &   1.69   &  4.22e-02 &   0.88  &  4.20e-04 &   1.83  &  4.22e-02 &   0.88  \\
 5 &     7.53e-04 &   1.73   &  2.27e-02 &   0.89  &  1.18e-04 &   1.83  &  2.27e-02 &   0.89  \\
 6 &     2.22e-04 &   1.76   &  1.21e-02 &   0.91  &  3.28e-05 &   1.84  &  1.21e-02 &   0.91  \\
 7 &     6.43e-05 &   1.79   &  6.42e-03 &   0.92  &  9.07e-06 &   1.85  &  6.42e-03 &   0.92  \\
 8 &     1.83e-05 &   1.81   &  3.38e-03 &   0.92  &  2.49e-06 &   1.86  &  3.38e-03 &   0.92  \\
 9 &     5.12e-06 &   1.83   &  1.77e-03 &   0.93  &  6.80e-07 &   1.87  &  1.77e-03 &   0.93  \\ \hline
\multicolumn{2}{l}{\text{Expected}} &   2     &           &   1     & &   2     &           &   1
 \end{array}
 \]
  \caption{Errors and experimental orders of convergence for the boundary value problem. In the bottom line, orders of convergence expected from \eqref{E4.14}.
  Graded mesh. $\mu = 0.66$}\label{T6.2}
\end{table}

\begin{table}
  \centering
 \[
  \begin{array}{c@{\;}|@{\;}c@{\;}c@{\;}|@{\;}c@{\;}c@{\;}||@{\;}c@{\;}c@{\;}|@{\;}c@{\;}c}
 j &    \Vert \bar y-\bar y_h(\bar u)\Vert_{L^2(\Omega)} &   EOC  &     \Vert \bar y-\bar y_h(\bar u)\Vert_{H^1(\Omega)}  &  EOC &    \Vert \bar \varphi-\tilde \varphi_h\Vert_{L^2(\Omega)} &   EOC  &     \Vert \bar \varphi-\tilde \varphi_h\Vert_{H^1(\Omega)}  &  EOC \\ \hline
 1 &     8.30e-02 &         &  2.48e-01 &          &  2.27e-02 &         &  2.65e-01 &        \\
 2 &     2.59e-02 &   1.68  &  1.40e-01 &   0.82   &  5.53e-03 &   2.03  &  1.42e-01 &   0.90 \\
 3 &     5.94e-03 &   2.12  &  7.04e-02 &   0.99   &  1.28e-03 &   2.12  &  7.05e-02 &   1.01 \\
 4 &     1.77e-03 &   1.75  &  3.74e-02 &   0.91   &  3.46e-04 &   1.88  &  3.74e-02 &   0.91 \\
 5 &     5.38e-04 &   1.72  &  2.00e-02 &   0.90   &  9.50e-05 &   1.86  &  2.00e-02 &   0.90 \\
 6 &     1.37e-04 &   1.98  &  1.01e-02 &   0.99   &  2.42e-05 &   1.97  &  1.01e-02 &   0.99 \\
 7 &     3.74e-05 &   1.87  &  5.22e-03 &   0.95   &  6.34e-06 &   1.93  &  5.22e-03 &   0.95 \\
 8 &     1.03e-05 &   1.85  &  2.71e-03 &   0.95   &  1.67e-06 &   1.92  &  2.71e-03 &   0.95 \\
 9 &     2.63e-06 &   1.98  &  1.36e-03 &   0.99   &  4.22e-07 &   1.99  &  1.36e-03 &   0.99 \\ \hline
\multicolumn{2}{l}{\text{Expected}} &   2     &           &   1  &  & 2     &           &   1
 \end{array}
 \]
  \caption{Errors and experimental orders of convergence for the boundary value problem. In the bottom line, orders of convergence expected from \eqref{E4.14}.
  Graded mesh. $\mu = 0.6$}\label{T6.3}
\end{table}

\begin{table}
  \centering
 \[
  \begin{array}{c@{\;}|@{\;}c@{\;}c@{\;}|@{\;}c@{\;}c@{\;}||@{\;}c@{\;}c@{\;}|@{\;}c@{\;}c}
 j &    \Vert \bar y-\bar y_h(\bar u)\Vert_{L^2(\Omega)} &   EOC  &     \Vert \bar y-\bar y_h(\bar u)\Vert_{H^1(\Omega)}  &  EOC &    \Vert \bar \varphi-\tilde \varphi_h\Vert_{L^2(\Omega)} &   EOC  &     \Vert \bar \varphi-\tilde \varphi_h\Vert_{H^1(\Omega)}  &  EOC \\ \hline
 1 &     8.30e-02 &         &  2.48e-01 &          &  2.27e-02 &         &  2.65e-01 &         \\
 2 &     1.80e-02 &   2.20  &  1.26e-01 &   0.98   &  4.50e-03 &   2.33  &  1.27e-01 &   1.06  \\
 3 &     4.20e-03 &   2.10  &  6.31e-02 &   0.99   &  1.10e-03 &   2.03  &  6.32e-02 &   1.00  \\
 4 &     1.07e-03 &   1.97  &  3.19e-02 &   0.99   &  2.83e-04 &   1.96  &  3.19e-02 &   0.99  \\
 5 &     2.75e-04 &   1.96  &  1.61e-02 &   0.99   &  7.22e-05 &   1.97  &  1.61e-02 &   0.99  \\
 6 &     7.03e-05 &   1.97  &  8.11e-03 &   0.99   &  1.83e-05 &   1.98  &  8.11e-03 &   0.99  \\
 7 &     1.78e-05 &   1.98  &  4.08e-03 &   0.99   &  4.60e-06 &   1.99  &  4.08e-03 &   0.99  \\
 8 &     4.51e-06 &   1.98  &  2.05e-03 &   0.99   &  1.16e-06 &   1.99  &  2.05e-03 &   0.99  \\
 9 &     1.14e-06 &   1.99  &  1.03e-03 &   1.00   &  2.90e-07 &   2.00  &  1.03e-03 &   1.00  \\ \hline
\multicolumn{2}{l}{\text{Expected}} &   2     &           &   1    &  & 2     &           &   1
 \end{array}
 \]
  \caption{Errors and experimental orders of convergence for the boundary value problem. In the bottom line, orders of convergence expected from \eqref{E4.14}.
  Graded mesh. $\mu = 0.5$}\label{T6.4}
\end{table}

Next, we turn to the control problem and check the estimate in Theorem \ref{T5.7}. Notice that we should obtain order of convergence $h$ for both graded-meshes and quasi-uniform meshes. We summarize the results in Table \ref{T6.5}.
\begin{table}
  \centering
  \[
  \begin{array}{c@{\;}|@{\;}c@{\;}c@{\;}|@{\;}c@{\;}c@{\;}||@{\;}c@{\;}c@{\;}|@{\;}c@{\;}c}
     \multicolumn{3}{c}{\text{Quasi-uniform}} &   \multicolumn{2}{c}{\mu=0.66} &   \multicolumn{2}{c}{\mu=0.6} &   \multicolumn{2}{c}{\mu=0.5}\\
 j &  \Vert \bar u-\bar u_h\Vert_{L^2(\Gamma)} &   EOC  &     \Vert \bar u-\bar u_h\Vert_{L^2(\Gamma)}  &  EOC  &     \Vert \bar u-\bar u_h\Vert_{L^2(\Gamma)}  &  EOC  &     \Vert \bar u-\bar u_h\Vert_{L^2(\Gamma)}  &  EOC \\ \hline
1 &  2.58e-01 &           &       2.58e-01 &          &   2.58e-01 &          &   2.58e-01 &          \\
2 &  1.33e-01 &     0.96  &       1.15e-01 &     1.16 &   1.15e-01 &     1.16 &   1.15e-01 &     1.17 \\
3 &  6.79e-02 &     0.97  &       6.05e-02 &     0.93 &   5.62e-02 &     1.03 &   5.61e-02 &     1.03 \\
4 &  3.45e-02 &     0.98  &       2.90e-02 &     1.06 &   2.90e-02 &     0.96 &   2.80e-02 &     1.01 \\
5 &  1.75e-02 &     0.98  &       1.48e-02 &     0.97 &   1.44e-02 &     1.00 &   1.41e-02 &     0.98 \\
6 &  8.80e-03 &     0.99  &       7.34e-03 &     1.01 &   7.22e-03 &     1.00 &   7.07e-03 &     1.00 \\
7 &  4.42e-03 &     0.99  &       3.68e-03 &     1.00 &   3.61e-03 &     1.00 &   3.53e-03 &     1.00 \\
8 &  2.22e-03 &     1.00  &       1.84e-03 &     1.00 &   1.80e-03 &     1.00 &   1.76e-03 &     1.00 \\
9 &  1.11e-03 &     1.00  &       9.19e-04 &     1.00 &   9.01e-04 &     1.00 &   8.82e-04 &     1.00 \\ \hline
\multicolumn{2}{l}{\text{Expected}}             &     1    &               &     1    &               &     1    &               &     1
\end{array}
 \]
  \caption{Errors and experimental orders of convergence for the optimal control problem. In the bottom line, orders of convergence expected from Theorem \ref{T5.7}.}\label{T6.5}
\end{table}

Note that in this example the regularity of the adjoint state is even $\bar\varphi\in W^{2,\infty}_{\vec\gamma}(\Gamma)$ for $\vec\gamma =(\gamma,0,0,0,0,0)$ with $\gamma >4/3$. This leads to superconvergence properties in the convergence in the norms of $L^2(\Omega)$ and $L^2(\Gamma)$ of both the state and adjoint state variable, where, despite expecting order of convergence $1$, as for the control, we obtain the same order of convergence as the one for the boundary value problem, i.e, $1.33$ or almost $2$ in our examples. This phenomenon will be studied in a future paper.


\end{document}